# INTEGRAL CURVES OF NOISY VECTOR FIELDS AND STATISTICAL PROBLEMS IN DIFFUSION TENSOR IMAGING: NONPARAMETRIC KERNEL ESTIMATION AND HYPOTHESES TESTING

By Vladimir Koltchinskii,[1] Lyudmila Sakhanenko[2] and Songhe Cai[3]

*The University of New Mexico, Michigan State University and The University of New Mexico*

Let $v$ be a vector field in a bounded open set $G \subset \mathbb{R}^d$. Suppose that $v$ is observed with a random noise at random points $X_i$, $i = 1, \ldots, n$, that are independent and uniformly distributed in $G$. The problem is to estimate the integral curve of the differential equation

$$\frac{dx(t)}{dt} = v(x(t)), \qquad t \geq 0, x(0) = x_0 \in G,$$

starting at a given point $x(0) = x_0 \in G$ and to develop statistical tests for the hypothesis that the integral curve reaches a specified set $\Gamma \subset G$. We develop an estimation procedure based on a Nadaraya–Watson type kernel regression estimator, show the asymptotic normality of the estimated integral curve and derive differential and integral equations for the mean and covariance function of the limit Gaussian process. This provides a method of tracking not only the integral curve, but also the covariance matrix of its estimate. We also study the asymptotic distribution of the squared minimal distance from the integral curve to a smooth enough surface $\Gamma \subset G$. Building upon this, we develop testing procedures for the hypothesis that the integral curve reaches $\Gamma$.

The problems of this nature are of interest in diffusion tensor imaging, a brain imaging technique based on measuring the diffusion tensor at discrete locations in the cerebral white matter, where the diffusion of water molecules is typically anisotropic. The diffusion tensor data is used to estimate the dominant orientations of the

Received June 2005; revised November 2006.
[1]Supported by NSF Grant DMS-03-04861 and NIH Grant NIBIB 1 R01 EB002618-01.
[2]Supported by MSU Grant IRGP 03-42179.
[3]Supported by NIH Grant NIBIB 1 R01 EB002618-01.
*AMS 2000 subject classifications.* 62G05, 62G10, 62G20, 62P10.
*Key words and phrases.* Diffusion tensor imaging, kernel regression estimator, vector fields, integral curves, functional central limit theorem.

diffusion and to track white matter fibers from the initial location following these orientations. Our approach brings more rigorous statistical tools to the analysis of this problem providing, in particular, hypothesis testing procedures that might be useful in the study of axonal connectivity of the white matter.

**1. Introduction.** Let $G \subset \mathbb{R}^d$ be a bounded open set. Suppose a vector field $v : G \mapsto \mathbb{R}^d$ is observed at points $X_i \in G, i = 1, \ldots, n$, with random errors; that is, the observations are

$$V_i = v(X_i) + \xi_i,$$

where $\xi, \xi_1, \ldots, \xi_n$ are i.i.d. bounded random vectors (r.v.) $\mathbb{E}\xi = 0$ and $\mathrm{Cov}(\xi, \xi) = \Sigma$.

We are interested in the Cauchy problem for the differential equation (ODE)

$$\frac{dx(t)}{dt} = v(x(t)), \qquad t \geq 0, x(0) = x_0 \in G, \tag{1.1}$$

which of course can be equivalently written in an integral form,

$$x(t) = x_0 + \int_0^t v(x(s))\, ds.$$

Our goal is to provide an estimate $\hat{X}(t), t \geq 0$, of its solution based on the data $(X_i, V_i), i = 1, \ldots, n$, and, most importantly, to study the asymptotic behavior as $n \to \infty$ of statistics such as $\inf_{0 \leq t \leq T} d^2(\hat{X}(t), \Gamma)$, where $\Gamma \subset G$ is a given subset of $G$ (most often, it will be the boundary of a specified region in $G$) and $d(x, \Gamma) := \inf\{|x - y| : y \in \Gamma\}$ is the usual Euclidean distance from $x$ to $\Gamma$. This would allow us to suggest tests of the hypothesis that the true trajectory $x(t), 0 \leq t \leq T$, reaches a certain region in $G$.

Our main interest in this problem is related to its potential applications to *diffusion tensor imaging* (DTI), a technique in brain research introduced several years ago and often combined with conventional MRI (see, e.g., [7]). The diffusion of water molecules at a given location is characterized by a symmetric positive definite $3 \times 3$ diffusion matrix (diffusion tensor). The principal eigenvector of this matrix shows the dominant direction of the diffusion. In cerebral white matter, the diffusion is typically anisotropic and DTI allows one to recover its dominant directions by measuring the diffusion tensor field within voxels at a discrete set of locations and computing principle eigenvectors of diffusion matrices (thus transforming the tensor field into a vector field, see Figure 1). The fiber tract then can be reconstructed by following the directions of the vectors in small steps from a specified initial location, which essentially means solving numerically the ODE generated by the vector field. This provides a noninvasive approach to study the axonal connectivity of white matter fiber inside a brain region. The method is often referred to as *white matter fiber tractography*.

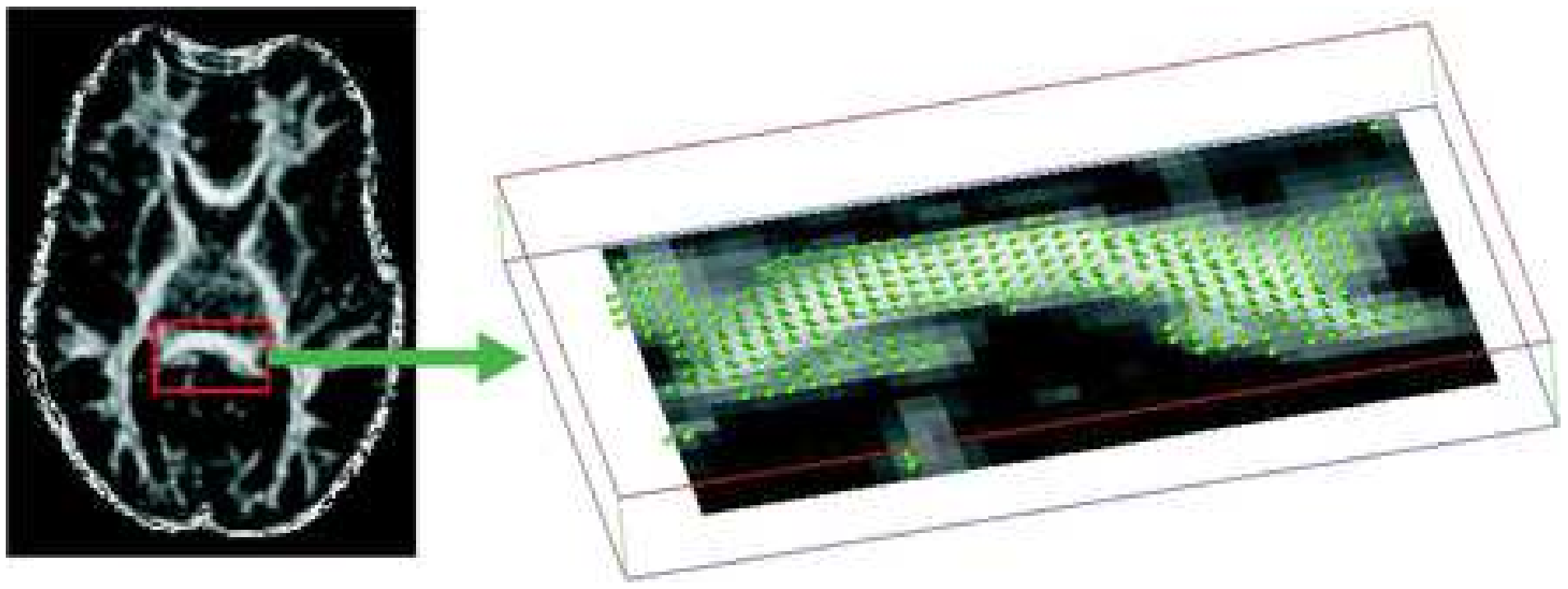

FIG. 1. *Shows the 3D vector field based on DTI data. The left graph is a fractional anisotropic (FA) map; the right graph gives the 3D visualization of the vector field inside a rectangular region of the FA map.*

Since the diffusion tensor field is being measured at a discrete set of locations and each matrix in the field represents an average within a voxel corrupted with noise, it becomes crucial to use some methods of smoothing of tensor or vector fields or regularization techniques that restrict fibers to smooth paths. For instance, Basser et al. [2] applied B-spline smoothing to the tensor field; Poupon et al. [18] used Markov random field models to obtain a regularized estimate of the vector field. However, even fiber track estimates involving smoothing would possess a certain degree of variability and very little is known about quantitative ways to assess the variability of fiber track estimates (although Parker, Barker and Buckley [14] and Jones [10] suggest some Monte Carlo and bootstrap approaches to the problem), which would facilitate the development of more rigorous approaches to statistical analysis of DTI data (the situation is somewhat different in conventional MRI and fMRI where approaches based on rather deep statistical understanding of the problem are becoming more common; see, e.g., [19] and [17]). This seems to be an important task especially because the mathematical models used in fiber tractography are rather involved and the existing methods utilize tools coming from very different areas (see [1, 2, 3, 6, 12, 13, 15, 16, 18, 20]).

The goal of the paper is to make the first (and rather modest) step toward better theoretical understanding of statistical problems in DTI. Our approach to estimation of $x(t), t \geq 0$, is based on the Nadaraya–Watson type kernel regression estimate (NWE) $\hat{V}(x), x \in G$, of the vector field $v(x), x \in G$, which then is substituted for $v$ into the ODE (1.1). The solution of the resulting Cauchy problem is an estimate $\hat{X}(t), t \geq 0$. (This approach is somewhat akin to that of Basser et al. [2]: the difference is that we are applying smoothing to the vector field and not to the tensor field; also, we are using kernel regression based smoothing instead of splines.) We establish in Section 2 (the proofs are given in Section 3) the asymptotic normality of this estimator, that is, the weak convergence (in the space of continuous functions) of the properly normalized deviation process $\hat{X}(t) - x(t), t \in [0, T]$, to a vector-valued Gaussian process on $[0, T]$ with mean and (matrix-valued) covariance function that depend on the vector field $v$, on the covariance matrix $\Sigma$ of the noise $\xi$ and on the kernel of NWE. We derive differential equations for the covariance function of the limit process, which allows us to develop a technique of *simultaneous tracking of the fiber path and its covariance* (see Section 4). We also study (Section 2) the asymptotic distribution (as $n \to \infty$) of the distance $\inf_{0 \le t \le T} d^2(\hat{X}(t), \Gamma)$ from the estimated integral curve $\hat{X}(t), t \in [0, T]$, to a set $\Gamma \subset G$. The asymptotic distributions for such distances happen to be especially simple in the case when the minimum of the function $[0, T] \ni t \mapsto d(x(t), \Gamma)$ is attained at a single point. In this case, the distributions are either normal or $\chi^2$-type, and they depend on the geometry of $\Gamma$ and on whether the true integral curve $x(t), t \in [0, T]$ reaches $\Gamma$ and in which way. These results allow one to bring into the analysis of DTI data some tools of rigorous statistical inference. In particular, one can use the asymptotic normality of $\hat{X}(t)$ to construct confidence ellipsoids for $x(t)$ for a fixed $t$; one can go further and try to use the results on Gaussian processes to develop nonparametric confidence bands and hypotheses tests for the whole integral curve $x(t), t \in [0, T]$; one can develop statistical tests for the hypothesis that the true integral curve $x(t), t \in [0, T]$, reaches a specified subregion of $G$; one can develop confidence intervals for the distance from $x(t), t \in [0, T]$, to a subregion. The last two possibilities are especially important since they are related to the problem of axonal connectivity which is one of the central issues in DTI. In Section 4, some of these possibilities are studied both for simulated and for real data.

There is a number of issues that (in our view) go beyond the scope of the paper, but that needs to be addressed to develop a comprehensive methodology based on our approach. First of all, the choice of NWE of the vector field $v$ is relatively arbitrary. Similar theory could, in principle, be developed for a number of other smoothing techniques. Moreover, it might be more natural and statistically more appealing to do smoothing of the underlying tensor field and only then to compute the principal eigenvectors creating a vector field. However, methods of perturbation theory will be needed to develop the asymptotic theory of such estimators, which would make the mathematical analysis of the problem more involved. Also, it is not common in DTI (at least, to our knowledge) to measure direction vectors at random locations, so, regression models with fixed design would be more appropriate than the model with random design (which is used in the paper primarily because the mathematics looks nicer). The vectors $V_i$ are usually unit eigenvectors of diffusion matrices, so, it would be more natural to consider nonparametric regression models for directional data rather than additive noise models. We are not exploring many important aspects of kernel-type nonparametric estimation (such as data-driven choice of the bandwidth, minimax lower bounds, optimal convergence rates, adaptation, etc.). Finally, in fiber tractography it is of great importance to take into account the possibility of fiber paths branching or intersecting one another. This is not covered by our model (because of the uniqueness of the solution of ODE) and the extension of our results to this case poses some nontrivial problems.

Realizing the importance of all these and some other issues, we, however, believe that the results we obtained so far might be of some interest for further development of a comprehensive statistical theory of DTI.

**2. A kernel estimate of integral curves and its asymptotic normality.** We will assume that $G \subset \mathbb{R}^d$ is a bounded open set of Lebesgue measure 1 and, for simplicity, that $X_i$ are i.i.d. uniformly distributed in $G$ and that the r.v.'s $\{\xi_i\}$ are independent of $\{X_i\}$. We also assume that $\operatorname{supp}(v) := \overline{\{x : v(x) \neq 0\}} \subset G$, which allows us to set $v = 0$ outside of $G$. Furthermore, we need a smoothness assumption on the vector field $v$. Unless stated otherwise, we assume that it is twice continuously differentiable.

We will use the following NWE of the vector field $v$ (see, e.g., Efromovich [4]):

$$\hat{V}(x) = \hat{V}_n(x) = \frac{1}{nh^d} \sum_{i=1}^{n} K\left(\frac{x - X_i}{h}\right) V_i,$$

with a kernel $K$ satisfying standard assumptions, in particular,

$$\int_{\mathbb{R}^d} K(x)\, dx = 1, \qquad \int_{\mathbb{R}^d} K(x) x\, dx = 0,$$

and with some bandwidth parameter $h = h_n$. It will be also convenient to assume that $K$ has bounded support, where it is twice continuously differentiable (the last assumption can be replaced by a more mild one in most of the results, but it is not of great importance in the context of the paper). As a result, the estimate $\hat{V}(x) = 0$ outside a bounded neighborhood of $G$. Comparing with the standard NWE, our estimate is simpler: since the distribution of $X_i$ is known (it is uniform), we do not need to use the kernel density estimator in the denominator of $\hat{V}$. In the case of nonuniform design some other smoothing techniques, such as local polynomial models (see [5]), might also be of interest.

Then, we define a plug-in estimate of the solution $x(t), t \geq 0$, as the solution $\hat{X}(t) = \hat{X}_n(t), t \geq 0$, of the Cauchy problem

$$\frac{d\hat{X}(t)}{dt} = \hat{V}(\hat{X}(t)), \qquad t \geq 0, \qquad \hat{X}(0) = x_0 \in G, \tag{2.1}$$

which is equivalent to the integral equation

$$\hat{X}(t) = x_0 + \int_0^t \hat{V}(\hat{X}(s))\,ds. \tag{2.2}$$

Note that since both $v$ and $\hat{V}$ vanish outside a neighborhood of $G$ ($v$ actually vanishes outside $G$ itself), $x(t)$ and $\hat{X}(t)$ will remain in this neighborhood for all $t > 0$.

To be specific, we assume that all vectors are vector columns; the sign $^*$ will denote transposition of vectors or matrices. Whenever it is convenient, we use the notation $\langle\cdot,\cdot\rangle$ for the inner product in $\mathbb{R}^d$ and $\mathbb{I}$ denotes the identity matrix. Also $\mathcal{N}(\mu,\Sigma)$ denotes the normal distribution with mean $\mu$ and covariance matrix $\Sigma$ and $Z \sim \mathcal{N}(\mu,\Sigma)$ denotes a r.v. with this distribution.

In what follows, we also need an estimate of the derivative of $v$ and we use for this purpose

$$\hat{V}'(x) = \frac{1}{nh^{d+1}} \sum_{i=1}^n V_i \biggl(K'\biggl(\frac{x - X_i}{h}\biggr)\biggr)^*.$$

Under the above assumptions, $\hat{V}, \hat{V}'$ are consistent estimates of $v, v'$ uniformly in $\mathbb{R}^d$ (see Lemma 1 below).

Our first goal is to show that under the assumptions $h \to 0$ and $nh^{d+3} \to \beta \ge 0$ the sequence of stochastic processes $\sqrt{nh^{d-1}}(\hat{X}(t) - x(t)), 0 \le t \le T$, converges weakly in the space $C[0,T] = C([0,T],\mathbb{R}^d)$ of $\mathbb{R}^d$-valued continuous functions on $[0,T]$ to the Gaussian process $\xi(t)$ satisfying the SDE

$$\begin{aligned} d\xi(t) &= \frac{\sqrt{\beta}}{2} \int_{\mathbb{R}^d} K(u)\langle v''(x(t))u, u\rangle\,du\,dt + v'(x(t))\xi(t)\,dt \\ &\quad + (\psi(v(x(t)))[\Sigma + v(x(t))v^*(x(t))])^{1/2}\,dW(t) \end{aligned} \tag{2.3}$$

with initial condition $\xi(0) = 0$, where $W(t), t \ge 0$, is a standard Brownian motion in $\mathbb{R}^d$,

$$\psi(v) := \int_{\mathbb{R}} \Psi(v\tau)\,d\tau, \qquad \Psi(y) := \int_{\mathbb{R}^d} K(z)K(z+y)\,dz.$$

Note also that $v''(x(t))$ involved in (2.3) is a $d \times d \times d$-tensor and $\langle v''(x(s))z, z\rangle$ is a vector-valued quadratic form. In what follows, $M_\beta(t)$ denotes the mean and $C(t)$ the covariance matrix of $\xi(t)$ (which does not depend on $\beta$). In Section 4 we provide differential equations for $M_\beta(t)$ and $C(t)$.

THEOREM 1. *Suppose that $h_n \to 0$ and $nh_n^{d+2} \to \infty$ as $n \to \infty$. Then for all $T > 0$*

$$\sup_{0 \le t \le T} |\hat{X}_n(t) - x(t)| \to 0 \qquad \text{as } n \to \infty,$$

*in probability. Suppose also that* $nh_n^{d+3} \to \beta \geq 0$ *as* $n \to \infty$. *Let* $T > 0$ *and suppose that for some* $\gamma = \gamma_T > 0$ *and for all* $0 \leq s \leq t \leq T$

$$\left| \frac{1}{t-s} \int_s^t v(x(\lambda))\, d\lambda \right| \geq \gamma.$$

*Then the sequence of stochastic processes* $\sqrt{nh_n^{d-1}}(\hat{X}_n(t) - x(t)), 0 \leq t \leq T$, *converges weakly in the space* $C[0,T]$ *to the Gaussian process* $\xi(t), 0 \leq t \leq T$.

REMARK. The condition of boundedness of the support of the kernel $K$ can be replaced by the conditions that the functions $\bar{\Psi}$ and $\Lambda$ (defined in the proof of Theorem 1) are integrable. The proof of Theorem 1 goes through and the theorem applies to such kernels as the Gaussian.

We turn now to some consequences of the asymptotic normality. In particular, we are interested in asymptotic properties of statistics of the type $\inf_{t\in[0,T]} m(d(\hat{X}(t), \Gamma))$, where $\Gamma$ is a subset of $G$, $d(x, \Gamma)$ is a distance from $x$ to $\Gamma$ and $m$ is a monotone function [e.g., $m(u) = u^2$ or $m(u) = u$, $u > 0$]. In other words, we want to study the asymptotic behavior of the minimal distance from the estimated integral curve $\hat{X}$ to a target set $\Gamma$. Such results are of statistical significance since they allow one to develop tests for hypotheses that the true integral curve $x(t)$ is passing through a given region or to construct confidence intervals for the distance to the region. We will study this problem under the assumption that the function $\varphi(x) := m(d(x, \Gamma))$ is smooth enough, which leads to a somewhat more general question about convergence in distribution (subject to a proper normalization) of the sequence $\inf_{t\in[0,T]} \varphi(\hat{X}(t)) - \inf_{t\in[0,T]} \varphi(x(t))$.

THEOREM 2. *Let* $x(t), t \geq 0$, *be an integral curve starting at* $x(0) = x_0 \in G$. *Suppose that* $\varphi : G \mapsto \mathbb{R}$ *is continuously differentiable. Denote*

$$M := \left\{ \tau \in [0,T] : \varphi(x(\tau)) = \inf_{0 \leq t \leq T} \varphi(x(t)) \right\}$$

*and suppose that* $M \subset (0,T)$. *Finally, suppose the conditions of Theorem 1 hold. Then the sequence of r.v.s*

$$\sqrt{nh_n^{d-1}} \left[ \inf_{t\in[0,T]} \varphi(\hat{X}(t)) - \inf_{t\in[0,T]} \varphi(x(t)) \right]$$

*converges in distribution to the r.v.* $\inf_{\tau\in M} \xi(\tau)^* \varphi'(x(\tau))$. *In particular, if the minimal set* $M$ *consists of only one point* $\tau \in (0,T)$, *then the above sequence is asymptotically normal with mean* $M_\beta(\tau)$ *and variance* $\sigma^2 = (\varphi'(x(\tau)))^* C(\tau) \varphi'(x(\tau))$. *Suppose now that* $\varphi$ *is twice continuously differentiable. If, for all* $\tau \in M$, $\varphi'(x(\tau)) = 0$ *and* $\varphi''(x(\tau))(v(x(\tau)), v(x(\tau))) > 0$, *then the sequence of r.v.s*

$$nh_n^{d-1} \left[ \inf_{t\in[0,T]} \varphi(\hat{X}(t)) - \inf_{t\in[0,T]} \varphi(x(t)) \right]$$

*converges in distribution to the r.v.*

$$\frac{1}{2}\inf_{\tau\in M}\biggl[\varphi''(x(\tau))(\xi(\tau),\xi(\tau)) - \frac{(\varphi''(x(\tau))(v(x(\tau)),\xi(\tau)))^2}{\varphi''(x(\tau))(v(x(\tau)),v(x(\tau)))}\biggr].$$

*If the minimal set consists only of one point $\tau$, then the limit becomes*

$$\frac{1}{2}\biggl[\varphi''(x(\tau))(Z,Z) - \frac{(\varphi''(x(\tau))(v(x(\tau)),Z))^2}{\varphi''(x(\tau))(v(x(\tau)),v(x(\tau)))}\biggr],$$

$$Z\sim\mathcal{N}(M_\beta(\tau),C(\tau)) \qquad \text{in } \mathbb{R}^d.$$

*On the other hand, if for all $u\in\mathbb{R}^d$, $\varphi''(x(\tau))(v(x(\tau)),u)=0$, then the distributional limit of the sequence*

$$nh_n^{d-1}\Bigl[\inf_{t\in[0,T]}\varphi(\hat{X}(t)) - \inf_{t\in[0,T]}\varphi(x(t))\Bigr]$$

*is $\frac{1}{2}\inf_{\tau\in M}\varphi''(x(\tau))(\xi(\tau),\xi(\tau))$, which in the unique minimum case is $\frac{1}{2}\varphi''(x(\tau))\times(Z,Z)$.*

Consider two typical examples.

Corollary 1. *Let $a\in G$ and $x(t), t\ge 0$, be an integral curve starting at $x(0)=x_0\in G$. Suppose that for some $\tau\in(0,T)$ $\inf_{0\le t\le T}|x(t)-a|^2=|x(\tau)-a|^2$, and, moreover, suppose that $\tau$ is the only point where the infimum is attained. Suppose also the conditions of Theorem 1 hold. If $x(\tau)\ne a$, then the sequence*

$$\sqrt{nh_n^{d-1}}\Bigl[\inf_{0\le t\le T}|\hat{X}(t)-a|^2 - \inf_{0\le t\le T}|x(t)-a|^2\Bigr]$$

*is asymptotically normal with mean $2M_\beta(\tau)^*(x(\tau)-a)$ and variance $\sigma^2=4(x(\tau)-a)^*C(\tau)(x(\tau)-a)$. If $x(\tau)=a$, then the sequence $nh_n^{d-1}\inf_{0\le t\le T}|\hat{X}(t)-a|^2$ converges in distribution to the r.v.*

$$|Z|^2-\frac{(v(x(\tau))^*Z)^2}{|v(x(\tau))|^2}, \qquad Z\sim\mathcal{N}(M_\beta(\tau),C(\tau)) \qquad \text{in } \mathbb{R}^d.$$

Corollary 2. *Let $\Gamma:=\{x:|x-a|=r\}\subset G$ be a sphere. Let*

$$d(x,\Gamma):=\inf_{y\in\Gamma}|x-y|=||x-a|-r|$$

*be the distance from $x$ to $\Gamma$ and let $x(t), t\ge 0$, be an integral curve starting at $x(0)=x_0\in G$. Suppose that for some $\tau\in(0,T)$*

$$\inf_{0\le t\le T}d^2(x(t),\Gamma)=d^2(x(\tau),\Gamma)=:D^2,$$

*and, moreover, suppose that $\tau$ is the only point where the infimum is attained. Suppose also the conditions of Theorem 1 hold. If $D^2 > 0$, then the sequence*

$$\sqrt{nh_n^{d-1}}\left[\inf_{0\le t\le T} d^2(\hat{X}(t),\Gamma) - D^2\right]$$

*is asymptotically normal with mean $2DM_\beta(\tau)^* n(x(\tau))$ and variance*

$$\sigma^2 = 4D^2 n(x(\tau))^* C(\tau) n(x(\tau)),$$

*where $n(x) := \frac{x-a}{|x-a|}$. If $D^2 = 0$ and, moreover, the vector $v(x(\tau))$ is tangent to $\Gamma$, then the sequence $nh_n^{d-1}\inf_{0\le t\le T} d^2(\hat{X}(t),\Gamma)$ converges in distribution to the r.v. $\gamma^2$, where $\gamma$ is a normal random variable with mean $M_\beta(\tau)^* n(x(\tau))$ and variance $n(x(\tau))^* C(\tau) n(x(\tau))$.*

REMARKS. 1. The result can be extended to more general smooth surfaces $\Gamma$. In this case, $n(x)$ would be the unit normal vector to $\Gamma$ at the point $x' \in \Gamma$ that is the closest to $x$ (assuming the uniqueness of such a point).

2. Suppose $H \subset G$ is an open nonempty subset of $G$ with boundary $\partial H = \Gamma$. Let $x(t), t \in [0,T]$, be the integral curve with initial condition $x(0) = x_0, x_0 \notin H \cup \Gamma$. If for some $t \in [0,T]$ $x(t) \in H$, then $\inf_{0\le t\le T} d^2(x(t),\Gamma) = 0$ since $x(t), t \in [0,T]$, is a continuous function. Also, it easily follows from the first statement of Theorem 1 that with probability tending to 1 we have $\inf_{0\le t\le T} d^2(\hat{X}(t),\Gamma) = 0$ (since $\hat{X}$, being close to $x$ uniformly in $[0,T]$, must enter the set $H$ and hence cross its boundary $\Gamma$). As a result, for *any* sequence $a_n \to \infty$

$$a_n\left[\inf_{0\le t\le T} d^2(\hat{X}(t),\Gamma) - \inf_{0\le t\le T} d^2(x(t),\Gamma)\right] = a_n \inf_{0\le t\le T} d^2(\hat{X}(t),\Gamma)$$

tends to 0 in probability and in distribution.

**3. Proofs of the main results.** We will need the following quite standard statement which is given without proof.

LEMMA 1. *Suppose that $h \to 0$ and $nh^{d+2} \to \infty$ as $n \to \infty$. Under the assumptions above, we have in probability*

$$\sup_{x\in\mathbb{R}^d} |\hat{V}(x) - \mathbb{E}\hat{V}(x)| \to 0,$$

$$\sup_{x\in\mathbb{R}^d} |\hat{V}(x) - v(x)| \to 0 \quad \text{and} \quad \sup_{x\in\mathbb{R}^d} |\hat{V}'(x) - v'(x)| \to 0.$$

PROOF OF THEOREM 1. *Consistency.* Let $y(t) := \hat{X}(t) - x(t)$. We have

$$y(t) = \int_0^t [\hat{V}(\hat{X}(s)) - v(x(s))]\,ds$$
$$= \int_0^t (\hat{V} - v)(\hat{X}(s))\,ds + \int_0^t [v(\hat{X}(s)) - v(x(s))]\,ds,$$

which implies (using a Lipschitz condition on $v$ and the fact that both $\hat{X}$ and $x$ remain in a bounded neighborhood of $G$) that with some constant $L$ for all $t \in [0, T]$

$$|y(t)| \le T \sup_{x \in \mathbb{R}^d} |\hat{V}(x) - v(x)| + L \int_0^t |y(s)|\,ds.$$

Using Lemma 1 and this Gronwall–Bellman inequality (see [8]), this easily implies consistency.

*Asymptotic representation of* $\hat{X} - x$. We will establish that

$$\hat{X}(t) - x(t) = z(t) + \delta(t), \qquad t \in [0, T], \tag{3.1}$$

where $z(t) = z_n(t), \delta(t) = \delta_n(t)$ are sequences of stochastic processes such that the sequence $\sqrt{nh^{d-1}} z_n(t), 0 \le t \le T$, converges in distribution to the Gaussian process $\xi(t)$ and

$$\sup_{0 \le t \le T} |\delta(t)| = o_p\left(\frac{1}{\sqrt{nh^{d-1}}}\right). \tag{3.2}$$

The following representation is obvious:

$$\begin{aligned} y(t) &= \int_0^t [\hat{V}(\hat{X}(s)) - v(x(s))]\,ds \\ &= \int_0^t (\hat{V} - v)(x(s))\,ds + \int_0^t v'(x(s)) \cdot y(s)\,ds + R(t), \end{aligned} \tag{3.3}$$

where the remainder is defined as

$$\begin{aligned} R(t) &:= \int_0^t [\hat{V}(\hat{X}(s)) - \hat{V}(x(s)) - v'(x(s)) \cdot y(s)]\,ds \\ &= \int_0^t [(\hat{V} - v)(\hat{X}(s)) - (\hat{V} - v)(x(s))]\,ds \\ &\quad + \int_0^t [v(\hat{X}(s)) - v(x(s)) - v'(x(s)) \cdot y(s)]\,ds. \end{aligned}$$

Note that

$$|(\hat{V} - v)(\hat{X}(s)) - (\hat{V} - v)(x(s))|$$

$$= \left| \int_0^1 (\hat{V} - v)'(a\hat{X}(s) + (1-a)x(s))\, da \cdot y(s) \right|$$
$$\leq \sup_{0 \leq a \leq 1} |(\hat{V} - v)'(a\hat{X}(s) + (1-a)x(s))| \cdot |y(s)|$$
$$\leq \sup_{x \in \mathbb{R}^d} |\hat{V}'(x) - v'(x)||y(s)|.$$

Also,

$$|v(\hat{X}(s)) - v(x(s)) - v'(x(s))y(s)|$$
$$= \left| \int_0^1 [v'(a\hat{X}(s) + (1-a)x(s)) - v'(x(s))]\, da \cdot y(s) \right|$$
$$\leq \sup_{0 \leq a \leq 1} |v'(a\hat{X}(s) + (1-a)x(s)) - v'(x(s))| \cdot |y(s)|$$
$$\leq r(|y(s)|) \cdot |y(s)|,$$

where $r(\delta) := \sup_{x \in \mathbb{R}^d} \sup_{|y| \leq \delta} |v'(x+y) - v'(x)| \to 0$ as $\delta \to 0$, since $v'$ is uniformly continuous on $G$. Then for all $t \in [0, T]$

$$|R(t)| \leq \left( \sup_{x \in \mathbb{R}^d} |\hat{V}'(x) - v'(x)| + r\left( \sup_{0 \leq s \leq T} |y(s)| \right) \right) \int_0^t |y(s)|\, ds. \tag{3.4}$$

Since $\sup_{0 \leq s \leq T} |y(s)| \to 0$ in probability and by Lemma 1 $\sup_{x \in \mathbb{R}^d} |\hat{V}'(x) - v'(x)| \to 0$ in probability, we have for all $T > 0$

$$\sup_{0 \leq t \leq T} |R(t)| = o_p\left( \int_0^T |y(s)|\, ds \right), \tag{3.5}$$

which also implies

$$\sup_{0 \leq t \leq T} |R(t)| = o_p\left( \sup_{0 \leq t \leq T} |y(t)| \right). \tag{3.6}$$

Denote by $z(t), t \geq 0$, the solution of the integral equation

$$z(t) := \int_0^t [\hat{V}(x(s)) - v(x(s))]\, ds + \int_0^t v'(x(s))z(s)\, ds.$$

In other words, $z$ satisfies the equation

$$\frac{dz(t)}{dt} = \hat{V}(x(t)) - v(x(t)) + v'(x(t))z(t), \qquad z(0) = 0.$$

Let $\delta(t) := y(t) - z(t)$. We have

$$\delta(t) = \int_0^t v'(x(s)) \cdot \delta(s)\, ds + R(t),$$

which implies

$$|\delta(t)| \le |R(t)| + \int_0^t |v'(x(s))| \cdot |\delta(s)|\, ds \le \sup_{0\le t\le T} |R(t)| + \int_0^t |v'(x(s))| \cdot |\delta(s)|\, ds.$$

Applying again the Gronwall–Bellman inequality, we get

$$|\delta(t)| \le \sup_{0\le t\le T} |R(t)| \exp\left\{\int_0^t |v'(x(u))|\, du\right\} \le C \sup_{0\le t\le T} |R(t)|, \qquad 0 \le t \le T,$$

with some constant $C > 0$, since the exponent is bounded. As a result, by (3.5),

$$\sup_{0\le t\le T} |\delta(t)| = o_p\left(\int_0^T |y(t)|\, dt\right) \qquad \text{as } n \to \infty.$$

Since $y(t) = z(t) + \delta(t)$, we also have

$$\sup_{0\le t\le T} |\delta(t)| = o_p\left(\int_0^T |z(t)|\, dt\right) \qquad \text{as } n \to \infty.$$

It follows from the weak convergence of $\sqrt{nh^{d-1}} z_n$ to the continuous stochastic process $\xi$ (established below) that

$$\int_0^T |z(t)|\, dt = O_p\left(\frac{1}{\sqrt{nh^{d-1}}}\right). \tag{3.7}$$

Therefore, we also have (3.2).

Consider now the following differential equation in $\mathbb{R}^d$:

$$\frac{du(t)}{dt} = \frac{dF(t)}{dt} + v'(x(t))u(t), \qquad u(0) = 0,$$

where $F$ is a continuously differentiable function from $[0,T]$ into $\mathbb{R}^d$ with $F(0) = 0$. This equation has the unique solution $u(t) = \mathcal{U}F(t)$, where $\mathcal{U}$ is a linear mapping from the space $C_0^{(1)}[0,T]$ of all continuously differentiable functions $F$ on $[0,T]$ with $F(0) = 0$ into $C[0,T]$. Another routine application of the Gronwall–Bellman inequality shows that $\mathcal{U}$ is a continuous (in fact, even Lipschitz) mapping with respect to the uniform distance.

Denote

$$\eta_n(t) := \sqrt{nh^{d-1}} \int_0^t (\hat{V}(x(s)) - v(x(s)))\, ds,$$
$$\xi_n(t) := \sqrt{nh^{d-1}} z_n(t), \qquad t \ge 0.$$

Then $\xi_n = \mathcal{U}\eta_n$. We will show that the sequence of stochastic processes $\eta_n$ converges weakly in the space $C[0,T]$ to the stochastic process $\eta$ satisfying

the SDE

$$d\eta(t) = \frac{\sqrt{\beta}}{2}\int K(u)\langle v''(x(t))u, u\rangle\, du\, dt$$
$$+ (\psi(v(x(t)))[\Sigma + v(x(t))v^*(x(t))])^{1/2}\, dW(t)$$

with initial condition $\eta(0) = 0$. (Here and in what follows we often do not write the limits of integration for the integrals over $\mathbb{R}$ or over $\mathbb{R}^d$.) Since $\mathcal{U}$ is a continuous mapping from $(C_0^{(1)}[0,T], \|\cdot\|_\infty)$ into $(C[0,T], \|\cdot\|_\infty)$, this implies that $\xi_n = \mathcal{U}\eta_n$ also converges weakly to the Gaussian process $\mathcal{U}\eta$, which satisfies (2.3) and the initial condition $\xi(0) = 0$. Hence, $\mathcal{U}\eta = \xi$.

*Computation of mean and covariance.* Define $f_t(s) := I_{[0,t]}(s)$. Then

$$\eta_n(t) = \sqrt{nh^{d-1}}\int f_t(s)\cdot[\hat{V}(x(s)) - v(x(s))]\, ds$$
$$= \sqrt{nh^{d-1}}\int\left(X_n(f_t) - \int f_t(s)\cdot v(x(s))\, ds\right),$$

where

$$X_n(f) := \int f(s)\cdot\hat{V}(x(s))\, ds$$
$$= \frac{1}{nh^d}\sum_{j=1}^n \int f(s)K\left(\frac{x(s) - X_j}{h}\right) ds\cdot(v(X_j) + \xi_j).$$

Let $L$ denote the set of all bounded functions $f$ on $\mathbb{R}$ such that the support of $f$ is a subset of $[0,T]$ and $f$ is continuous almost everywhere in $\mathbb{R}$. Note that $L$ is a linear space and $f_t \in L, t \in [0,T]$. In asymptotic representations for the expectation and covariance of $X_n(f)$, we assume that $f \in L$. We start with $\mathbb{E}X_n(f)$:

$$\mathbb{E}X_n(f) = \frac{1}{h^d}\int f(s)\mathbb{E}K\left(\frac{x(s) - X}{h}\right)\cdot v(X)\, ds$$
$$= \frac{1}{h^d}\int f(s)\int K\left(\frac{x(s) - y}{h}\right)\cdot v(y)\, dy\, ds$$
$$= \int f(s)\int K(z)\cdot v(x(s) - zh)\, dz\, ds$$
$$(3.8)\qquad = \int f(s)\int K(z)\Big[v(x(s)) - hv'(x(s))\cdot z$$
$$+ \frac{h^2}{2}\langle v''(x(s))z, z\rangle + o(h^2)\Big]\, dz\, ds$$

$$= \int f(s)\cdot v(x(s))\,ds - h\int f(s)\cdot v'(x(s))\cdot\left(\int K(z)z\,dz\right)ds$$
$$+ \frac{h^2}{2}\int f(s)\cdot\int K(z)\langle v''(x(s))z,z\rangle\,dz\,ds + o(h^2),$$

where we have used the substitution $z=\frac{x(s)-y}{h}$. Note that under the assumption $nh^{d+3}\to\beta\ge 0$ we have

$$\mathbb{E}\eta_n(t)\to\frac{\sqrt{\beta}}{2}\int_0^t\int K(z)\langle v''(x(s))z,z\rangle\,dz\,ds=\mathbb{E}\eta(t).$$

Also,

$$\operatorname{Cov}(X_n(f),X_n(g))=\frac{1}{nh^{2d}}\operatorname{Cov}\left(\int f(s)K\left(\frac{x(s)-X}{h}\right)ds\cdot(v(x)+\xi),\right.$$
$$\left.\int g(s)K\left(\frac{x(s)-X}{h}\right)ds\cdot(v(X)+\xi)\right).$$

Since

$$\operatorname{Cov}\left(\int f(s)K\left(\frac{x(s)-X}{h}\right)ds\cdot v(X),\int g(s)K\left(\frac{x(s)-X}{h}\right)ds\cdot\xi\right)=0,$$
$$\operatorname{Cov}\left(\int f(s)K\left(\frac{x(s)-X}{h}\right)ds\cdot\xi,\int g(s)K\left(\frac{x(s)-X}{h}\right)ds\cdot v(X)\right)=0,$$

we have

$$\begin{aligned}
&\operatorname{Cov}(X_n(f),X_n(g))\\
&\quad=(\mathrm{I})+(\mathrm{II})\\
&\quad=\frac{1}{nh^{2d}}\operatorname{Cov}\left(\int f(s)K\left(\frac{x(s)-X}{h}\right)ds\cdot\xi,\int g(s)K\left(\frac{x(s)-X}{h}\right)ds\cdot\xi\right)\\
&\qquad+\frac{1}{nh^{2d}}\operatorname{Cov}\left(\int f(s)K\left(\frac{x(s)-X}{h}\right)ds\cdot v(X),\right.\\
&\qquad\qquad\left.\int g(s)K\left(\frac{x(s)-X}{h}\right)ds\cdot v(X)\right).
\end{aligned}$$

To handle (I), we write

$$\begin{aligned}
(\mathrm{I})&=\frac{1}{nh^{2d}}\mathbb{E}\left\{\int f(s)K\left(\frac{x(s)-X}{h}\right)ds\cdot\xi\cdot\xi^*\int K\left(\frac{x(u)-X}{h}\right)g(u)\,du\right\}\\
&=\frac{1}{nh^{2d}}\int\int\mathbb{E}\left\{K\left(\frac{x(s)-X}{h}\right)K\left(\frac{x(u)-X}{h}\right)\right\}f(s)\Sigma g(u)\,ds\,du.
\end{aligned}$$

Note that

$$\begin{aligned}&\mathbb{E}\Big\{K\Big(\frac{x(s)-X}{h}\Big)K\Big(\frac{x(u)-X}{h}\Big)\Big\}\\ &\qquad=\int K\Big(\frac{x(s)-y}{h}\Big)K\Big(\frac{x(u)-y}{h}\Big)\,dy\\ &\qquad=h^d\int K(z)K\Big(z+\frac{x(u)-x(s)}{h}\Big)\,dz\\ &\qquad=h^d\Psi\Big(\frac{x(u)-x(s)}{h}\Big).\end{aligned}$$

Changing variable as $u=s+\tau h$, we get

$$\text{(I)}=\frac{1}{nh^{d-1}}\int\int\Psi\Big(\frac{x(s+\tau h)-x(s)}{h}\Big)f(s)\Sigma g(s+\tau h)\,d\tau\,ds. \tag{3.9}$$

Note that $(x(s+\tau h)-x(s))/h\to v(x(s))$ as $n\to\infty$ and also for all $\tau$ and a.s. for $s$ $g(s+\tau h)\to g(s)$ as $n\to\infty$ (recall that the functions $f,g\in L$ and hence are continuous a.e. in $\mathbb{R}$). By assumptions $K$ has bounded support, implying that the support of $\Psi$ is also bounded. At the same time, we have

$$0<\gamma\le\Big|\frac{1}{u-s}\int_s^u v(x(\lambda))\,d\lambda\Big|\le\sup_{x\in\mathbb{R}^d}|v(x)|<+\infty.$$

Therefore, the function

$$\tau\mapsto\bar\Psi(\tau)=\sup_{0\le s\le u\le T}\Psi\Big(\tau\frac{1}{u-s}\int_s^u v(x(\lambda))\,d\lambda\Big)$$

also has bounded support and, since it is bounded, it is integrable. Thus, we can use Lebesgue dominated convergence to prove that

$$\begin{aligned}&\int\int\Psi\Big(\frac{x(s+\tau h)-x(s)}{h}\Big)f(s)\Sigma g(s+\tau h)\,d\tau\,ds\\ &\qquad\to\int\int\Psi(v(x(s))\tau)\,d\tau\,f(s)\Sigma g(s)\,ds=\int\psi(v(x(s)))f(s)\Sigma g(s)\,ds,\end{aligned}$$

which along with (3.9) yields

$$\text{(I)}=\frac{1+o(1)}{nh^{d-1}}\int\psi(v(x(s)))f(s)\Sigma g(s)\,ds.$$

[Indeed, the integration with respect to $s$ is in a finite range, the function $(s,\tau)\mapsto f(s)\Sigma g(s+\tau h)$ is uniformly bounded and

$$\Big|\Psi\Big(\frac{x(s+\tau h)-x(s)}{h}\Big)\Big|\le\bar\Psi(\tau),$$

so the dominated convergence can be used under the assumption that $\bar{\Psi}$ is integrable in $\mathbb{R}$.]

Similarly, the expression (II) can be written as

$$\begin{aligned}(\mathrm{II}) = \frac{1}{nh^{2d}} \int\int f(s)\mathbb{E}\Big\{K\Big(\frac{x(s)-X}{h}\Big)K\Big(\frac{x(u)-X}{h}\Big)v(X)\cdot v^*(X)\Big\}g(u)\,ds\,du \\ -\frac{1}{n}\int f(s)\cdot v(x(s))\,ds\int g(u)\cdot v(x(u))\,du(1+o(1)).\end{aligned}$$

Note that

$$\begin{aligned}&\mathbb{E}\Big\{K\Big(\frac{x(s)-X}{h}\Big)K\Big(\frac{x(u)-X}{h}\Big)v(X)\cdot v^*(X)\Big\} \\ &\qquad = \int K\Big(\frac{x(s)-y}{h}\Big)K\Big(\frac{x(u)-y}{h}\Big)v(y)\cdot v^*(y)\,dy \\ &\qquad = h^d\int K(z)K\Big(z+\frac{x(u)-x(s)}{h}\Big)v(x(s)-zh)\cdot v^*(x(s)-zh)\,dz,\end{aligned}$$

where we use the substitution $z=\frac{x(s)-y}{h}, dy = h^d dz$. Therefore,

$$\begin{aligned}&\frac{1}{nh^{2d}}\mathbb{E}\Big\{\int f(s)K\Big(\frac{x(s)-X}{h}\Big)\,ds\cdot v(X)\cdot v^*(X)\int K\Big(\frac{x(u)-X}{h}\Big)g(u)\,du\Big\} \\ &\qquad = \frac{1}{nh^d}\int\int f(s)\int K(z)K\Big(z+\frac{x(u)-x(s)}{h}\Big)v(x(s)-zh) \\ &\qquad\qquad \times v^*(x(s)-zh)\,dz\,g(u)\,ds\,du,\end{aligned}$$

which after the change of variable $u = s+\tau h$ becomes

$$\begin{aligned}\frac{1}{nh^{d-1}}\int\int f(s)\int K(z)K\Big(z+\frac{x(s+\tau h)-x(s)}{h}\Big)v(x(s)-zh) \\ \times v^*(x(s)-zh)\,dz\,g(s+\tau h)\,d\tau\,ds.\end{aligned}$$

As before, we use the dominated convergence (under the same conditions) to show that the last expression is equal to

$$\begin{aligned}&\frac{1+o(1)}{nh^{d-1}}\int\int f(s)\int K(z)K(z+\tau v(x(s)))\,dz\,d\tau\,v(x(s))\cdot v^*(x(s))g(s)\,ds \\ &\qquad = \frac{1+o(1)}{nh^{d-1}}\int f(s)\psi(v(x(s)))v(x(s))\cdot v^*(x(s))g(s)\,ds,\end{aligned}$$

implying that

$$(\mathrm{II}) = \frac{1+o(1)}{nh^{d-1}}\int f(s)\psi(v(x(s)))v(x(s))\cdot v^*(x(s))g(s)\,ds.$$

Finally, the covariance

$$
\begin{aligned}
&\mathrm{Cov}(X_n(f), X_n(g)) \\
&= (\mathrm{I}) + (\mathrm{II}) \qquad (3.10)\\
&= \frac{1+o(1)}{nh^{d-1}} \int \psi(v(x(s))) f(s) \cdot [\Sigma + v(x(s)) \cdot v^*(x(s))] \cdot g(s)\, ds.
\end{aligned}
$$

Thus,

$$
\begin{aligned}
\mathrm{Cov}(\eta_n(t_1), \eta_n(t_2)) &\to \mathrm{Cov}(\eta(t_1), \eta(t_2)) \\
&= \int_0^{t_1 \wedge t_2} \psi(v(x(s)))[\Sigma + v(x(s)) \cdot v^*(x(s))]\, ds.
\end{aligned}
$$

*Convergence of finite dimensional distributions* (*f.d.d.*). We now turn to the proof of asymptotic normality of $\eta_n(t), 0 \le t \le T$, in the sense of convergence of finite dimensional distributions. First we show that for all $f \in L$

$$
\sqrt{nh^{d-1}}\Big(X_n(f) - \int f(s) v(x(s))\, ds\Big)
$$

converges to a normal distribution. Since by (3.8)

$$
\begin{aligned}
&\sqrt{nh^{d-1}}\Big(\mathbb{E}X_n(f) - \int f(s) v(x(s))\, ds\Big) \\
&\quad \to \frac{\sqrt{\beta}}{2} \int f(s) \cdot \int K(z) \langle v''(x(s)) z, z\rangle\, dz\, ds,
\end{aligned}
$$

it is enough to establish the CLT for

$$
\sqrt{nh^{d-1}}(X_n(f) - \mathbb{E}X_n(f)) = \frac{1}{\sqrt{nh^{d+1}}} \sum_{j=1}^n (\chi_j - \mathbb{E}\chi_j),
$$

where

$$
\chi_j := \int f(s) K\Big(\frac{x(s) - X_j}{h}\Big)\, ds \cdot (v(X_j) + \xi_j).
$$

Under the assumptions we have made it is easy to check Lyapunov's conditions for the CLT, and to this end we bound the fourth moment of $\chi_j$,

$$
\begin{aligned}
\mathbb{E}|\chi_j|^4 &= \mathbb{E}(\chi_j^* \chi_j)^2 \\
&= \mathbb{E}\Big(\int\int K\Big(\frac{x(s) - X_j}{h}\Big) K\Big(\frac{x(s_1) - X_j}{h}\Big) f(s) \\
&\qquad \times (v(X_j) + \xi_j)^* (v(X_j) + \xi_j) f(s_1)\, ds\, ds_1\Big)^2.
\end{aligned}
$$

Under the assumption that $v$ and $\xi$ are bounded, this gives with some constant $C > 0$,

$$\begin{aligned}
\mathbb{E}|\chi_j|^4 &\le C\mathbb{E}\biggl(\int\int K\biggl(\frac{x(s)-X}{h}\biggr)K\biggl(\frac{x(s_1)-X}{h}\biggr)|f(s)||f(s_1)|\,ds\,ds_1\biggr)^2 \\
&= C\mathbb{E}\int_{\mathbb{R}^4} K\biggl(\frac{x(s)-X}{h}\biggr)K\biggl(\frac{x(s_1)-X}{h}\biggr)K\biggl(\frac{x(s_2)-X}{h}\biggr)K\biggl(\frac{x(s_3)-X}{h}\biggr) \\
&\qquad\qquad \times |f(s)||f(s_1)||f(s_2)||f(s_3)|\,ds\,ds_1\,ds_2\,ds_3.
\end{aligned}$$

By change of variable, we then get

$$\begin{aligned}
\mathbb{E}|\chi_j|^4 &\le Ch^d\int_{\mathbb{R}^5} K(z)K\biggl(z+\frac{x(s_1)-x(s)}{h}\biggr)K\biggl(z+\frac{x(s_2)-x(s)}{h}\biggr) \\
&\qquad \times K\biggl(z+\frac{x(s_3)-x(s)}{h}\biggr)\,dz\,|f(s)||f(s_1)||f(s_2)| \\
&\qquad \times |f(s_3)|\,ds\,ds_1\,ds_2\,ds_3 \\
&= Ch^{d+3}\int_{\mathbb{R}^5} K(z)K\biggl(z+\tau_1\frac{x(s+\tau_1 h)-x(s)}{\tau_1 h}\biggr) \\
&\qquad \times K\biggl(z+\tau_2\frac{x(s+\tau_2 h)-x(s)}{\tau_2 h}\biggr) \\
&\qquad \times K\biggl(z+\tau_3\frac{x(s+\tau_3 h)-x(s)}{\tau_3 h}\biggr) \\
&\qquad \times dz|f(s)||f(s+\tau_1 h)||f(s+\tau_2 h)| \\
&\qquad \times |f(s+\tau_3 h)|\,ds\,d\tau_1\,d\tau_2\,d\tau_3.
\end{aligned}$$

Denote

$$\begin{aligned}
\Lambda(\tau_1,\tau_2,\tau_3) := \sup\int K(z)K\biggl(z+\tau_1\frac{x(s_1)-x(s)}{s_1-s}\biggr) \\
\times K\biggl(z+\tau_2\frac{x(s_2)-x(s)}{s_2-s}\biggr)K\biggl(z+\tau_3\frac{x(s_3)-x(s)}{s_3-s}\biggr)\,dz,
\end{aligned}$$

where the supremum is taken over all $s, s_1, s_2, s_3 \in [0,T]$. It follows from the conditions that the function $\Lambda$ is integrable in $\mathbb{R}^3$.

As a result, we get

$$\begin{aligned}
\mathbb{E}|\chi_j|^4 \le Ch^{d+3}\int\int\int \Lambda(\tau_1,\tau_2,\tau_3)\biggl(\int |f(s)|^4\,ds\biggr)^{1/4}\biggl(\int |f(s+\tau_1 h)|^4\,ds\biggr)^{1/4} \\
\times\biggl(\int |f(s+\tau_2 h)|^4\,ds\biggr)^{1/4}
\end{aligned}$$

$$\times \left( \int |f(s+\tau_3 h)|^4 \, ds \right)^{1/4} d\tau_1 \, d\tau_2 \, d\tau_3$$

$$= Ch^{d+3} \int |f(s)|^4 \, ds \iiint \Lambda(\tau_1, \tau_2, \tau_3) \, d\tau_1 \, d\tau_2 \, d\tau_3.$$

It follows that with some constant $C$

$$\frac{1}{n^2 h^{2(d+1)}} \sum_{j=1}^{n} \mathbb{E}|\chi_j - \mathbb{E}\chi_j|^4 \leq \frac{Cnh^{d+3}}{n^2 h^{2(d+1)}} = \frac{C}{nh^{d-1}} \to 0,$$

implying Lyapunov's conditions for the CLT. This shows the asymptotic normality of

$$\sqrt{nh^{d-1}}(X_n(f) - \mathbb{E}X_n(f)) \quad \text{and} \quad \sqrt{nh^{d-1}}\left( X_n(f) - \int f(s) v(x(s)) \, ds \right)$$

for all $f \in L$. Hence, if $f_1, \ldots, f_m \in L$ (which is a linear space), the CLT holds for any linear combination of $f_1, \ldots, f_m$. Using the standard characteristic function argument, this shows that the joint distribution of $(X_n(f_1), \ldots, X_n(f_m))$ is also asymptotically normal. Applying this to $f = f_t$ proves the convergence of finite dimensional distributions (f.d.d.) of the stochastic processes $\eta_n(t), 0 \leq t \leq T$, to f.d.d. of the Gaussian process $\eta(t), 0 \leq t \leq T$.

*Asymptotic equicontinuity.* To prove the weak convergence of the sequence of processes $\eta_n(t), 0 \leq t \leq T$, in the functional space $C[0,T]$, it remains to check the asymptotic equicontinuity condition. Since

$$\eta_n(t) = \sqrt{nh^{d-1}}(X_n(f_t) - \mathbb{E}X_n(f_t)) + \sqrt{nh^{d-1}}\left( \mathbb{E}X_n(f_t) - \int f_t(s) v(x(s)) \, ds \right)$$

and the bias term $\sqrt{nh^{d-1}}(\mathbb{E}X_n(f_t) - \int f_t(s) v(x(s)) \, ds)$ tends to $\mathbb{E}\eta(t)$ uniformly in $t \in [0,T]$ due to (3.8), we have to consider only the process $\zeta_n(t) := \sqrt{nh^{d-1}}(X_n(f_t) - \mathbb{E}X_n(f_t))$. To this end, we bound the fourth moment of $X_n(f) - \mathbb{E}X_n(f)$,

$$\begin{aligned} &\mathbb{E}|X_n(f) - \mathbb{E}X_n(f)|^4 \\ &\quad = \frac{1}{n^4 h^{4d}} \mathbb{E}\left| \sum_{j=1}^{n} (\chi_j - \mathbb{E}\chi_j) \right|^4 \\ &\quad = \frac{1}{n^4 h^{4d}} \left[ \frac{n(n-1)}{2} (\mathbb{E}|\chi - \mathbb{E}\chi|^2)^2 + n\mathbb{E}|\chi - \mathbb{E}\chi|^4 \right]. \end{aligned} \tag{3.11}$$

As before,

$$\begin{aligned} &\mathbb{E}|\chi - \mathbb{E}\chi|^2 \\ &\quad \leq \mathbb{E}|\chi|^2 \end{aligned}$$

$$
\begin{aligned}
&= \mathbb{E}\int\int K\biggl(\frac{x(s)-X}{h}\biggr)K\biggl(\frac{x(s_1)-X}{h}\biggr)f(s)\\
&\qquad\times (v(X)+\xi)^*(v(X)+\xi)f(s_1)\,ds\,ds_1\\
&\le C\int\int \mathbb{E}K\biggl(\frac{x(s)-X}{h}\biggr)K\biggl(\frac{x(s_1)-X}{h}\biggr)|f(s)||f(s_1)|\,ds\,ds_1\\
&\le Ch^d\int\int\int K(z)K\biggl(z+\frac{x(s_1)-x(s)}{h}\biggr)\,dz\,|f(s)||f(s_1)|\,ds\,ds_1\\
&\le Ch^{d+1}\int\int\int K(z)K\biggl(z+\tau\frac{x(s+\tau h)-x(s)}{\tau h}\biggr)\,dz\\
&\qquad\qquad\times |f(s)||f(s+\tau h)|\,ds\,d\tau\\
&\le Ch^{d+1}\int\int \bar{\Psi}(\tau)|f(s)||f(s+\tau h)|\,ds\,d\tau\\
&\le Ch^{d+1}\int \bar{\Psi}(\tau)\biggl(\int |f(s)|^2\,ds\biggr)^{1/2}\biggl(\int |f(s+\tau h)|^2\,ds\biggr)^{1/2}\,d\tau\\
&\le Ch^{d+1}\int \bar{\Psi}(\tau)\,d\tau\int |f(s)|^2\,ds.
\end{aligned}
$$

Plugging the bounds on $\mathbb{E}|\chi-\mathbb{E}\chi|^2$ and on $\mathbb{E}|\chi-\mathbb{E}\chi|^4$ in (3.11) yields with a large enough constant $C>0$,

$$
\begin{aligned}
&\mathbb{E}(\sqrt{nh^{d-1}}|X_n(f)-\mathbb{E}X_n(f)|)^4\\
&\qquad\le C\biggl[\frac{n^2h^{2d+2}}{n^2h^{2d+2}}\biggl(\int |f(s)|^2\,ds\biggr)^2+\frac{nh^{d+3}}{n^2h^{2d+2}}\int |f(s)|^4\,ds\biggr]\\
&\qquad\le C\biggl[\biggl(\int |f(s)|^2\,ds\biggr)^2+\frac{1}{nh^{d-1}}\int |f(s)|^4\,ds\biggr].
\end{aligned}
$$

We will apply it to $f:=f_{t_1}-f_{t_2}$ with $t_1,t_2\in[0,T]$. It easily follows that with some $L>0$, $\int |f_{t_1}(s)-f_{t_2}(s)|^2\,ds\le L|t_1-t_2|$ and $\int |f_{t_1}(s)-f_{t_2}(s)|^4\,ds\le L|t_1-t_2|$. Therefore we have (with some $C>0$)

$$
\mathbb{E}|\zeta_n(t_1)-\zeta_n(t_2)|^4\le C\biggl[|t_1-t_2|^2+\frac{1}{nh^{d-1}}|t_1-t_2|\biggr],
$$

which gives $\mathbb{E}|\zeta_n(t_1)-\zeta_n(t_2)|^4\le 2C|t_1-t_2|^2$ for $|t_1-t_2|\le\frac{1}{nh^{d-1}}$. If now $A_n$ is a maximal $\frac{1}{nh^{d-1}}$-separated subset of $[0,T]$, then the standard Kolmogorov type of chaining argument shows that for all $\varepsilon>0$

$$
\lim_{\delta\to 0}\limsup_{n\to\infty}\mathbb{P}\Bigl\{\sup_{t_1,t_2\in A_n,|t_1-t_2|\le\delta}|\zeta_n(t_1)-\zeta_n(t_2)|\ge\varepsilon\Bigr\}=0. \tag{3.12}
$$

Let $\pi_n$ be a mapping from $[0,T]$ into $A_n$ such that $\forall t\in[0,T]:|t-\pi_n t|\le 1/nh^{d-1}$. Using the definition of $\zeta_n(t)$, we easily get (with some constant

$C > 0$)

$$|\zeta_n(t_1) - \zeta_n(t_2)| \le CT\sqrt{nh^{d-1}} \sup_{x\in\mathbb{R}^d} |\hat{V}(x) - \mathbb{E}\hat{V}(x)||t_1 - t_2|.$$

Therefore,

$$\sup_{t\in[0,T]} |\zeta_n(t) - \zeta_n(\pi_n t)| \le CT\frac{1}{\sqrt{nh^{d-1}}} \sup_{x\in\mathbb{R}^d} |\hat{V}(x) - \mathbb{E}\hat{V}(x)|.$$

Using Lemma 1,

$$\sup_{t\in[0,T]} |\zeta_n(t) - \zeta_n(\pi_n t)| = o_P(1). \tag{3.13}$$

It immediately follows from (3.12) and (3.13) that

$$\lim_{\delta\to 0} \limsup_{n\to\infty} \mathbb{P}\left\{ \sup_{t_1,t_2\in[0,T]|t_1-t_2|\le\delta} |\zeta_n(t_1) - \zeta_n(t_2)| \ge \varepsilon \right\} = 0,$$

which is the asymptotic equicontinuity condition for the process $\zeta_n$. □

The proof of Theorem 2 can be found in [11].

## 4. Numerical implementation and examples.

4.1. *Remarks on numerical implementation.* It is not hard to show that the mean of the limiting Gaussian process $\xi(t)$ defined by (2.3) can be written as $M_\beta(t) = \sqrt{\beta}M(t)$, where $M$ satisfies the differential equation

$$\frac{dM(t)}{dt} = v'(x(t))M(t) + \frac{1}{2}\int K(z)\langle v''(x(t))z, z\rangle\, dz, \qquad M(0) = 0. \tag{4.1}$$

The covariance matrix $C(t)$ of $\xi(t)$ does not depend on $\beta$ and it satisfies the ODE

$$\begin{aligned} \frac{dC(t)}{dt} &= \psi(v(x(t)))[\Sigma + v(x(t)) \cdot v^*(x(t))] \\ &\quad + v'(x(t))C(t) + C(t)v'(x(t))^*, \\ C(0) &= 0. \end{aligned} \tag{4.2}$$

It is also easy to derive partial differential equations for the covariance function of the process $\xi(t)$, but they are not used in what follows.

We will use a simple Euler type method to solve the ODEs numerically (obviously, more sophisticated numerical methods can also be useful here, with potential improvement of the results). Let $\delta$ be a step size. Then the following recurrent relationships provide an approximation of equations (4.1) and (4.2):

$$\hat{X}_{k+1} := \hat{X}_k + \hat{V}(\hat{X}_k)\delta \qquad \text{with } \hat{X}_0 = x_0, \tag{4.3}$$

$$\begin{aligned} \hat{C}_0 &:= 0, \\ \hat{C}_{k+1} &:= \hat{C}_k + \delta[\psi(\hat{V}(\hat{X}_k))(\hat{\Sigma} + \hat{V}(\hat{X}_k)\hat{V}(\hat{X}_k)^*) \\ &\qquad + \hat{V}'(\hat{X}_k)\hat{C}_k + \hat{C}_k\hat{V}'(\hat{X}_k)^*], \end{aligned} \tag{4.4}$$

where $\hat{\Sigma} := n^{-1}\sum_{j=1}^n (V_i - \hat{V}(X_i))(V_i - \hat{V}(X_i))^*$ is an estimate of the covariance matrix $\Sigma$ of the noise $\xi_i$. Consistency of the estimator $\hat{\Sigma}$ easily follows from Lemma 1. In practice, the noise $\xi_i$ is not necessarily homogeneous and it might make sense to use localized versions of the above estimate. Obviously, the recurrent relationships (4.3) and (4.4) can be solved simultaneously, so, in fact, our approach is based on simultaneous tracking of the "fiber path" and its covariance matrix. We are doing this for $k = 1, \dots, N$, $N := [\frac{T}{\delta}]$.

It easily follows from the definition of the function $\psi$ (see Section 2) that if the kernel $K$ is spherically symmetric, then $\psi$ is also spherically symmetric, so it is a constant on the unit sphere in $\mathbb{R}^d$. In applications, the vector field $v$ consists of unit vectors. Hence, for a spherically symmetric kernel $K$, the $\psi$-factor in the ODE for $C(t)$ and in (4.4) can be replaced by a constant, simplifying the equations. In what follows, we use the standard Gaussian kernel $K$, which of course is spherically symmetric.

To estimate the function $M(t)$, one needs an estimator of $v''$. This can be done, for instance, by utilizing kernel estimators one more time. The entries of the estimator $\hat{V}''$ (which is a $d \times d \times d$-tensor) are defined as

$$\hat{V}''_{jkl}(x) = \frac{1}{n\tilde{h}^{d+2}} \sum_{i=1}^n \frac{\partial^2 K}{\partial x_j\,\partial x_k}\left(\frac{x - X_i}{\tilde{h}}\right) V_i^{(l)},$$

$V_i^{(l)}, l = 1, \dots, d$, being the components of the vector $V_i$. The kernel $K$ can be taken to be the same as in the estimate of $\hat{V}$, but the bandwidth parameter $\tilde{h} = \tilde{h}_n$ is different (so $\hat{V}''$ is not the second derivative of $\hat{V}$). To make $\hat{V}''$ a consistent estimator of $v''$ the assumptions $\tilde{h} \to 0$ and $n\tilde{h}^{d+4} \to \infty$ are needed. The second assumption does not hold for the bandwidth $h$ needed in Theorem 1. If $K$ is the standard Gaussian kernel, then

$$\hat{W}(x) := \int K(z)\langle \hat{V}''(x)z, z\rangle\,dz = \frac{1}{n\tilde{h}^{d+2}} \sum_{i=1}^n \left(\left|\frac{x - X_i}{\tilde{h}}\right|^2 + d\right) K\left(\frac{x - X_i}{\tilde{h}}\right) V_i.$$

The following recurrence relation [i.e. to be solved simultaneously with (4.3) and (4.4)] provides a numerical approximation of equation (4.1):

$$\hat{M}_{k+1} = \hat{M}_k + \delta[\hat{V}'(\hat{X}_k)\hat{M}_k + \tfrac{1}{2}\hat{W}(\hat{X}_k)] \qquad \text{with } \hat{M}_0 = 0. \tag{4.5}$$

Solving (4.3), (4.4) and (4.5) yields numerical approximations of $\hat{X}(t)$, $\hat{C}(t)$, and $\hat{M}(t), 0 \le t \le T$, that can be now used to compute $\min_{1\le k\le N} d^2(\hat{X}_k, \Gamma)$, which is a numerical approximation of $\inf_{0\le t\le T} d^2(\hat{X}(t), \Gamma)$, for a given set $\Gamma$, and also to compute other quantities needed for implementation of testing procedures. If the above minimum is attained at $\hat{k}$ and $\hat{\tau} := \hat{k}\delta$, then $\hat{\tau}$ can be used as an estimate of $\tau$ for which the minimal distance from the true integral curve $x(t), 0 \le t \le T$, to $\Gamma$ is attained. If such a $\tau$ is unique (as was assumed in Corollaries 1 and 2), then it is not hard to show consistency of $\hat{\tau}$ (under proper assumptions on $\delta$). This provides an approximation of the limit distributions in Corollaries 1 and 2.

The above considerations allow us to implement the testing procedures based on Corollaries 1 and 2. For instance, in the case of Corollary 1, the test statistic is approximated by

$$\hat{\Lambda} := nh^{d-1} \min_{1\le k\le N} |\hat{X}_k - a|^2. \tag{4.6}$$

Given a significance level $\alpha \in (0,1)$, the hypothesis that the integral curve $x(t), 0 \le t \le T$, passes through the point $a$ (against the alternative that it does not) is rejected if $\hat{\Lambda} \ge \Lambda_\alpha$, where $\Lambda_\alpha$ is determined from the equation $\mathbb{P}\{\bar{\Lambda} \ge \Lambda_\alpha\} = \alpha$. Here

$$\bar{\Lambda} := |Z|^2 - \frac{(\hat{V}(\hat{X}_{\hat{k}})^* Z)^2}{|\hat{V}(\hat{X}_{\hat{k}})|^2}, \qquad Z \sim \mathcal{N}(\sqrt{\beta}\hat{M}_{\hat{k}}, \hat{C}_{\hat{k}}) \qquad \text{in } \mathbb{R}^d.$$

We assume that $h = (\frac{\beta}{n})^{1/(d+3)}$ with $\beta > 0$; we set $\beta = 0$ if $h = h_n$ is such that $nh_n^{d+3} \to 0$.

Corollaries 1 and 2 can be also used to derive asymptotic approximations of the power of the test and to study how it depends on $D$ (the minimal distance from the true integral curve to $\Gamma$). For instance, in the case of Corollary 1, the power can be approximated by the expression

$$1 - \Phi\biggl(\frac{(nh^{d-1})^{-1/2}\Lambda_\alpha - (nh^{d-1})^{1/2}D^2 - 2\sqrt{\beta}DM(\tau)^* n(x(\tau))}{2D(n(x(\tau))^* C(\tau) n(x(\tau)))^{1/2}}\biggr),$$

where $\Phi$ is the standard normal distribution function, $D^2 := \inf_{0\le t\le T} |x(t) - a|^2$ and $n(x) := \frac{x-a}{|x-a|}$. Replacing $M(\tau)$, $C(\tau)$ and $x(\tau)$ by their "estimates" leads to the following expression describing the dependence of the power on the true distance $D$:

$$1 - \Phi\biggl(\frac{(nh^{d-1})^{-1/2}\Lambda_\alpha - (nh^{d-1})^{1/2}D^2 - 2\sqrt{\beta}D\hat{M}_{\hat{k}}^* n(\hat{X}_{\hat{k}})}{2D(n(\hat{X}_{\hat{k}})^* \hat{C}_{\hat{k}} n(\hat{X}_{\hat{k}}))^{1/2}}\biggr). \tag{4.7}$$

We are not addressing in any detail an important problem of choosing the bandwidth parameter $h$. For a fixed $t$ and $h = (\frac{\beta}{n})^{1/(d+3)}$ the asymptotic formula for the mean squared error matrix of $\hat{X}$ is (see Theorem 1)

$$\mathbb{E}(\hat{X}(t) - x(t))(\hat{X}(t) - x(t))^* \\ \approx n^{-4/(d+3)}[C(t)\beta^{-4/(d+3)} + M(t)M(t)^* \beta^{4/(d+3)}]$$

(note that the convergence rate $n^{-4/(d+3)}$ is optimal in a minimax sense provided that the vector field $v$ is twice continuously differentiable; this can be shown, e.g., following the approach of [9], Theorem IV.5.1). This immediately implies the following formula for the mean integrated squared error:

$$\mathbb{E}\int_0^T |\hat{X}(t) - x(t)|^2\,dt \\ \approx n^{-4/(d+3)}\biggl[\int_0^T \mathrm{Tr}(C(t))\,dt\,\beta^{-(d-1)/(d+3)} \\ + \int_0^T \mathrm{Tr}(M(t)M(t)^*)\,dt\,\beta^{4/(d+3)}\biggr],$$

which can be easily minimized with respect to $\beta$, and the minimal point $\bar{\beta}$ can be estimated based on the data (using the estimates of $C$ and $M$). Since one might be interested in optimizing not the global deviation of $\hat{X}$ from $x$ but rather the distance from $x$ to a set $\Gamma$ (as in Corollaries 1 and 2), an alternative is to use the asymptotics of these corollaries rather than the global result of Theorem 1. For instance, based on Corollary 1, the following asymptotic formula might be used:

$$\mathbb{E}\Bigl[\inf_{0\le t\le T}|\hat{X}(t) - a|^2 - \inf_{0\le t\le T}|x(t) - a|^2\Bigr]^2 \\ \approx n^{-4/(d+3)}[4(x(\tau) - a)^* C(\tau)(x(\tau) - a)\beta^{-(d-1)/(d+3)} \\ + 4(M(\tau)^*(x(\tau) - a))^2 \beta^{4/(d+3)}],$$

which again can be easily minimized with respect to $\beta$, and the minimal point $\bar{\beta}_1$ can be estimated based on the data. One can also try to develop an approach based on maximizing the power of the hypothesis tests considered above.

4.2. *Several experiments with simulated and real data.* We turn now to some of the results of our experiments with simulated and real data. First, we simulated a vector field with circular integral curves (Figure 2). It was observed at a finite number of random points uniformly distributed inside a rectangular domain in $\mathbb{R}^2$ with random noise. We used NWE to smooth the vector field and then computed an estimate of an integral curve starting at a given point by solving numerically the ODE generated by the smoothed field using Euler's method. Simultaneously with tracking the estimate of the integral curve, we have also tracked the covariance matrix of the estimate and used it to plot the 95 percent confidence ellipses along the integral curve. The results are shown in Figure 2 (see [11] for more experiments of similar nature).

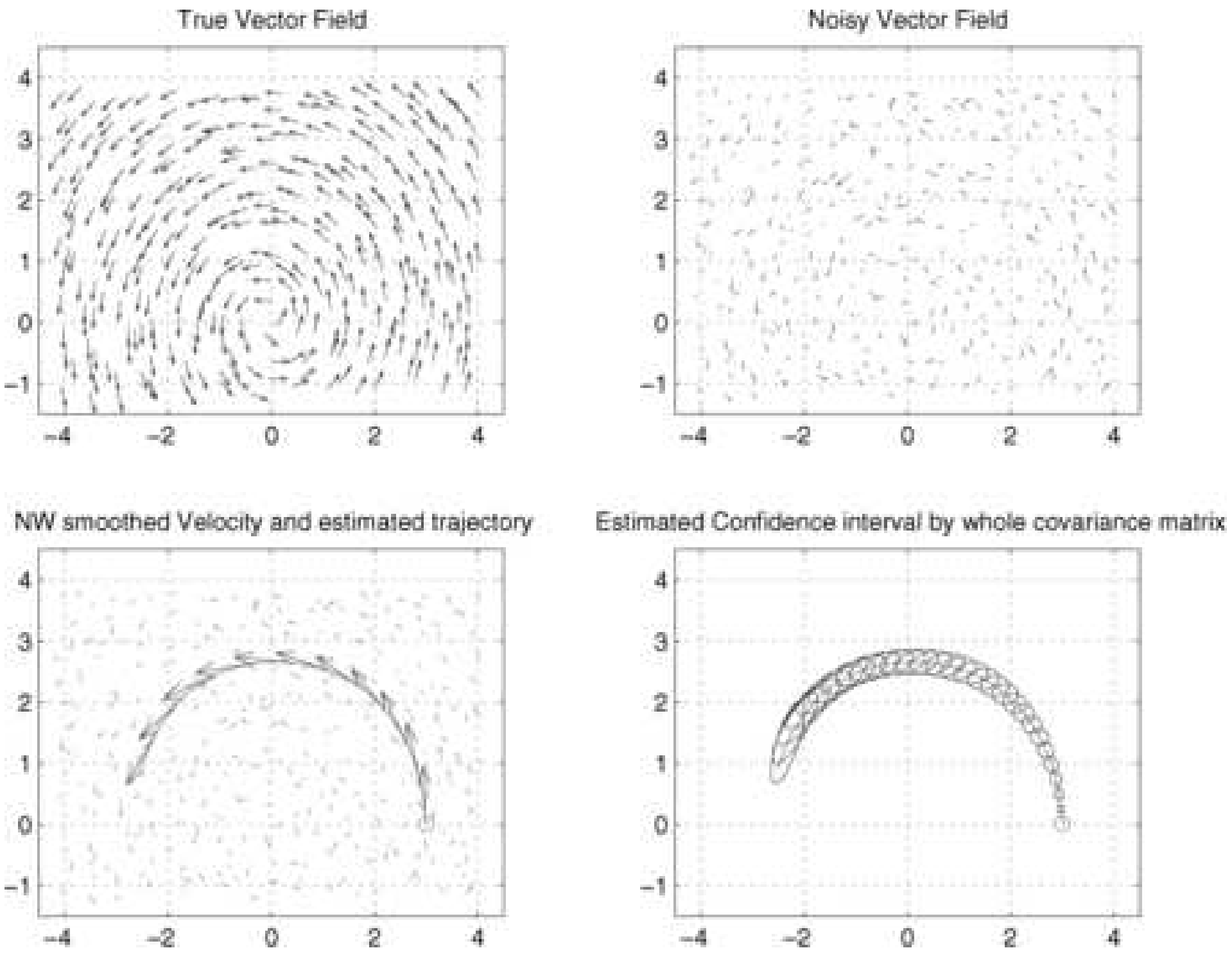

FIG. 2. *The top left figure shows the true vector field whose $x$- and $y$-coordinates are given by formulas $v_x = -\frac{y}{\sqrt{x^2+y^2}}, v_y = \frac{x}{\sqrt{x^2+y^2}}$ ("the circular field"). The top right figure represents the noisy vector field obtained by adding to the true field independent copies of $0.5Z,\ Z \sim \mathcal{N}(0, \mathbb{I})$ in $\mathbb{R}^2$. The bottom left figure shows the results of smoothing of the noisy vector field and integral curve estimation using NWE. The total number of points in the rectangle $n = 322$. The NWE was computed with $h = 0.85$ and the step size used in the numerical solution of ODE was $\delta = 0.02$. Finally, the bottom right figure depicts 95 percent confidence ellipses along the integral curve.*

Next, by means of Monte Carlo simulation we studied the accuracy of normal approximation of the distribution of the distance from the estimated integral curve to a given point or to a given sphere (see Corollaries 1, 2). To this end, we simulated the random points and the noisy vector field as in Figure 2 and computed the estimated integral curve based on NW regression smoothing. We repeated these simulations independently $N = 2000$ times and each time computed the square of the distance $\hat{D}^2$ to the point

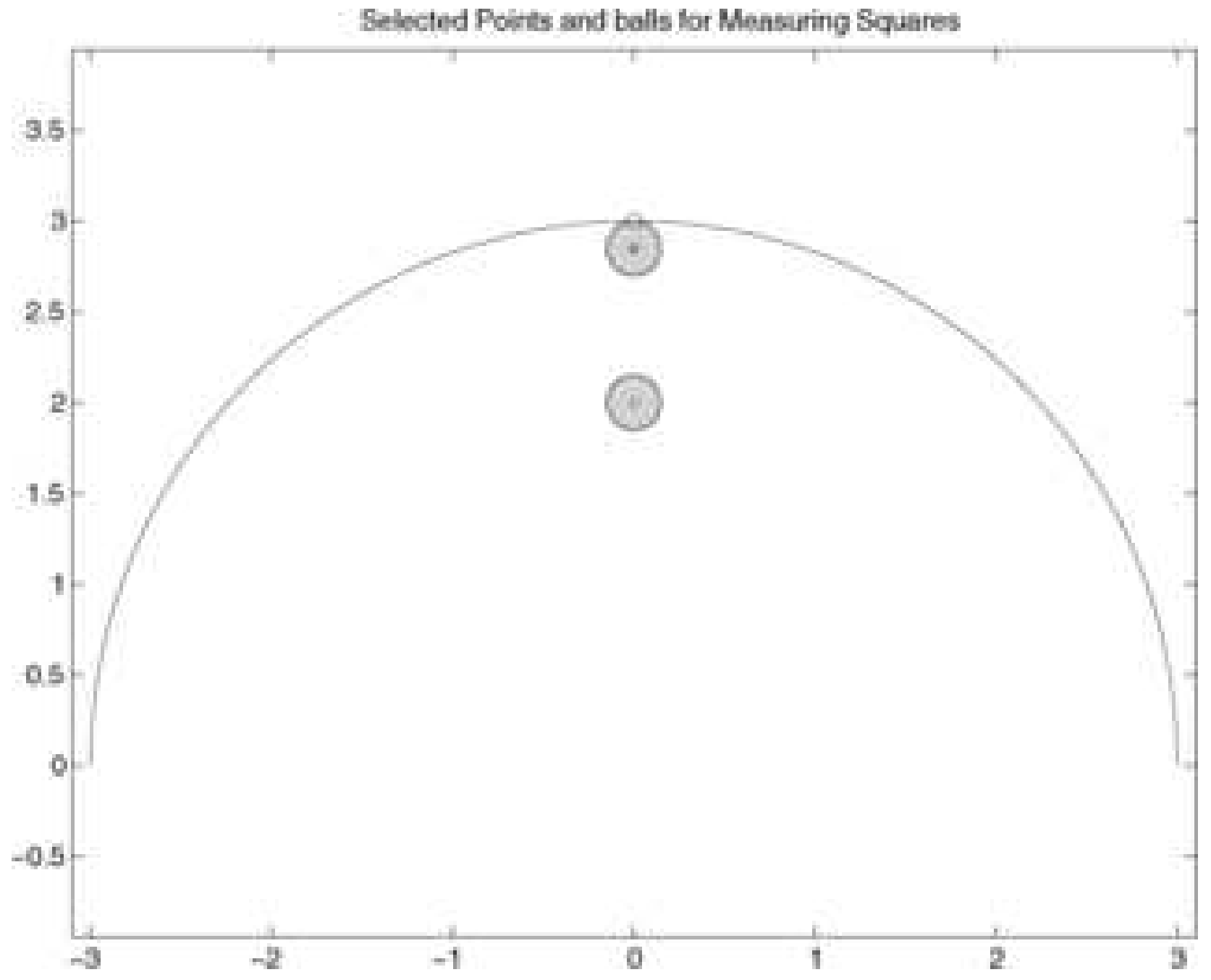

Fig. 3. *Shows a circular true integral curve and locations of points and balls of interest in a Monte Carlo study of the distribution of distances (see Figures 4 and 5 below).*

$(0,2)$ (labeled with $+$ in Figure 3). The squared distance from this point to the true integral curve was $D^2 = 1$. Each time we also computed the estimate $\hat{\sigma}^2$ of the variance $\sigma^2$ (see Corollary 1) and the standardized version of $\hat{D}^2$, given by the expression $\frac{\sqrt{nh}}{\hat{\sigma}}(\hat{D}^2 - D^2)$ (recall that $d = 2$ in our case). The histogram of the last variable is shown in the top part of Figure 4 in comparison with the standard normal curve. The bottom part of Figure 4 shows the results of a similar simulation in the case of the distance from the estimated integral curve to a sphere (a circle in our case; see Corollary 2). There is deviation of the histograms from normality that is quite understandable for a number of reasons: the fact that we ignored the bias $M_\beta$ in the normal approximation; in the case when $D^2 = 0$, Corollary 1 suggests that the asymptotic distribution should be of $\chi^2$-type rather than normal and because of this for small values of $D^2$ one can start seeing some deviations from normality for a finite sample; the variance $\sigma^2$ needed in normalization was replaced by its estimate $\hat{\sigma}^2$; the numerical approximation we are using to compute the distance has a certain impact on the distribution; and, last but not least, the sample size $n$ in our simulations is rather small for this type of CLT ($n = 77$). The Kolmogorov–Smirnov test clearly shows that these deviations from normality are very significant ($p \le 0.001$). However, when $n = 500$ the $p$-value of the test becomes of the order 0.0542 and for larger values of $n$ the deviations from normality are no longer statistically significant.

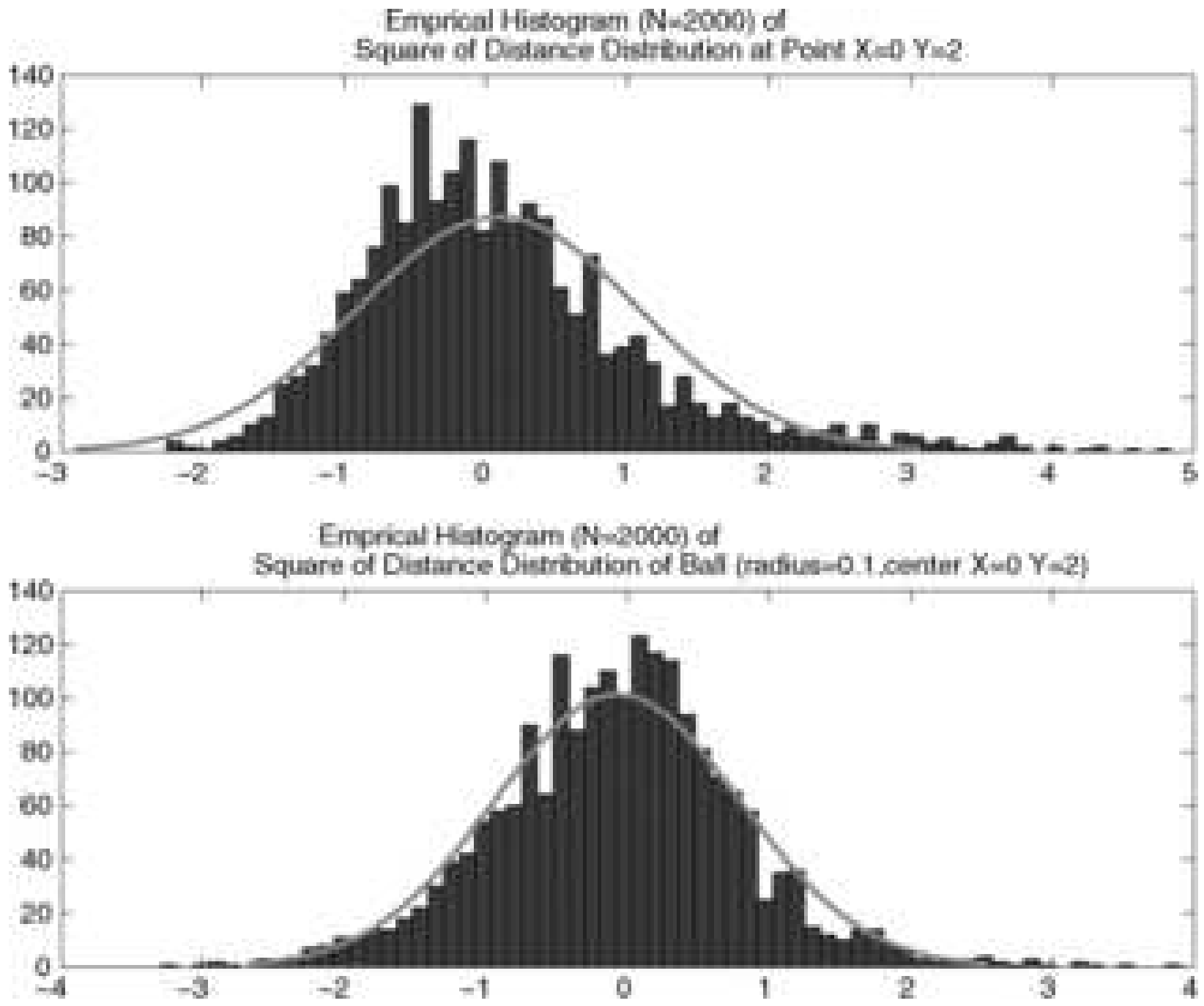

FIG. 4. *The top figure presents the histogram of standardized minimal squared distances from the estimated integral curves to the point $x = (0, 2)$ obtained by the Monte Carlo simulations ($N = 2000$); the bottom figure shows the histogram of standardized minimal squared distances from the estimated integral curve to the ball with center $x = (0, 2)$ and radius 0.1 obtained again by the Monte Carlo simulations ($N = 2000$). The histograms are compared with the standard normal distribution.*

Quite similarly, Figure 5 shows histograms of squared distances from the estimated integral curve to a specified point (the top figure) or to a specified circle (the bottom figure) in the case when the true integral curve passes through the point or is tangent to the circle. In this case, according to Corollaries 1 and 2, the asymptotic distribution of the squared distance should be of $\chi^2$-type.

Next we studied the power of testing the null hypothesis that the integral curve passes through a specified point of interest. The test is based on the second statement of Corollary 1. The test statistic is $\hat{\Lambda}$ given by (4.6). The top part of Figure 6 shows the true integral curve and also ten points of interest: one of them is on the curve (so that the null hypothesis is satisfied for this point) and nine other points represent alternatives. We estimated this integral curve based on $n = 77$ observations of a noisy vector field as in

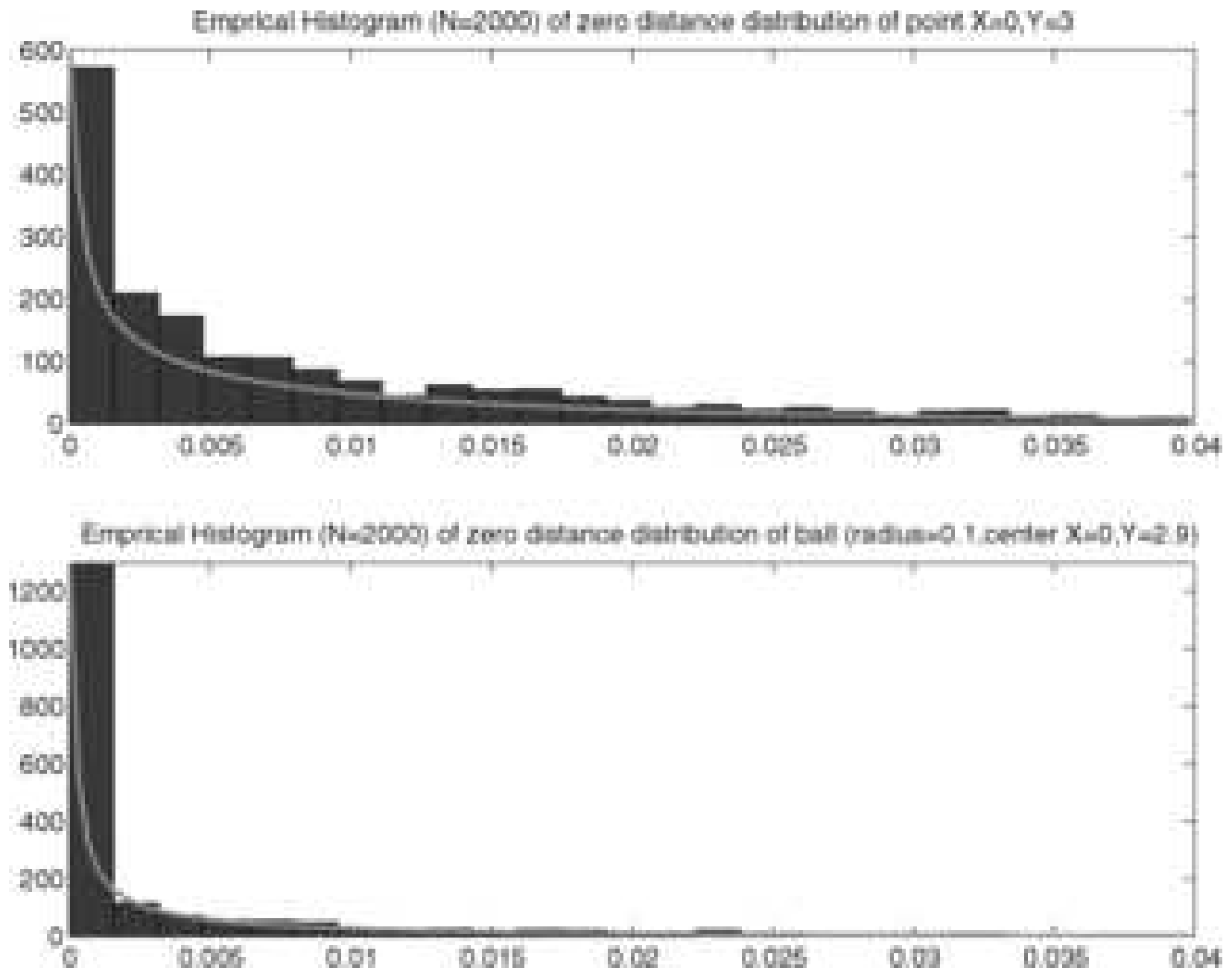

FIG. 5. *The top figure presents the histogram of the minimal squared distances from the estimated integral curve to the point* $x = (0, 3)$ *obtained by the Monte Carlo simulations* $(N = 2000)$ *in comparison with a* $\chi^2$*-type curve based on the theory. Note that now the point is on the true integral curve. The bottom figure shows the histogram of the minimal squared distances from the estimated integral curve to the ball with center* $x = (0, 2.9)$ *and radius 0.1 obtained by the Monte Carlo simulations* $(N = 2000)$ *again in comparison with a* $\chi^2$*-type curve based on the theory. The empirical distributions in this case are much closer to* $\chi^2$*-type than the distributions shown in Figure 4.*

Figure 2. We repeated the experiment 1000 times, each time simulating the data, estimating the integral curve and testing the hypothesis with significance level $\alpha = 0.05$. The red curve shown in the bottom part of Figure 6 represents the empirical estimation of the power of our test (the frequency of rejecting the null hypothesis) for each of the alternatives. The blue curve represents the value of the power based on the theoretical formula (4.7) (which seems to consistently overestimate the power).

Figures 7(a)–(c) give some examples of fiber tracking and visualization of 95% confidence ellipsoids for real DTI data and Figure 8 represents what we call *the p-value map*. It can be used to assess the degree of connectivity of points with a given path.

Although it is not the goal of this paper to develop a comprehensive methodology of statistical analysis of DTI data that takes into account more complicated issues, such as crossings or branchings of fibers, some of our results suggest possible ways to address these problems (as well as the problem of stopping criteria of the tracking). Usually the vector fields involved in DTI problems satisfy the condition $|v(x)| = 1$ since this is the field of unit principal eigenvectors of diffusion matrices. If the vector field is smooth around a given point $x \in G$, then $\hat{V}(x)$, as a weighted local average of unit vectors $v(X_i)$ of approximately the same direction plus a small noise, should have norm close to 1 (for large enough $n$). On the other hand, if two or more trajectories intersect at the point $x$, then $\hat{V}(x)$ becomes a weighted local average of unit vectors of several different directions and, as a result, $|\hat{V}(x)|$ will be, with high probability, significantly smaller than 1. Based on these heuristics, one can try to test the null hypothesis that a trajectory $x(t)$ starting at a given point $x_0$ *does not intersect* another trajectory at a given moment $t > 0$. It is natural to use the deviation $|\hat{V}(\hat{X}(t))|^2 - 1$ as a test statistic and to reject the null hypothesis if the value of this statistic is small enough. To make it more precise, one can use the following asymptotic result that easily follows from the CLT via a version of the delta method (and the result of Theorem 1, is also needed to handle the remainder terms): under the assumptions of Theorem 1, if $h_n \to 0$, $nh_n^{d+2} \to \infty$ and $nh_n^{d+3} \to 0$ as $n \to \infty$, then the sequence of r.v.s $\sqrt{nh_n^d}(|\hat{V}(\hat{X}(t))|^2 - 1)$ converges in distribution to the normal random variable with mean 0 and variance $\sigma^2 := 4\int K^2(u)\,du(1+v(x(t))^*\Sigma v(x(t)))$, which can be replaced by a plug-in estimate $\hat{\sigma}^2$ in a straightforward way. The statistic

$$\nu_n := \sqrt{nh_n^d}\frac{|\hat{V}(\hat{X}(t))|^2 - 1}{\hat{\sigma}}$$

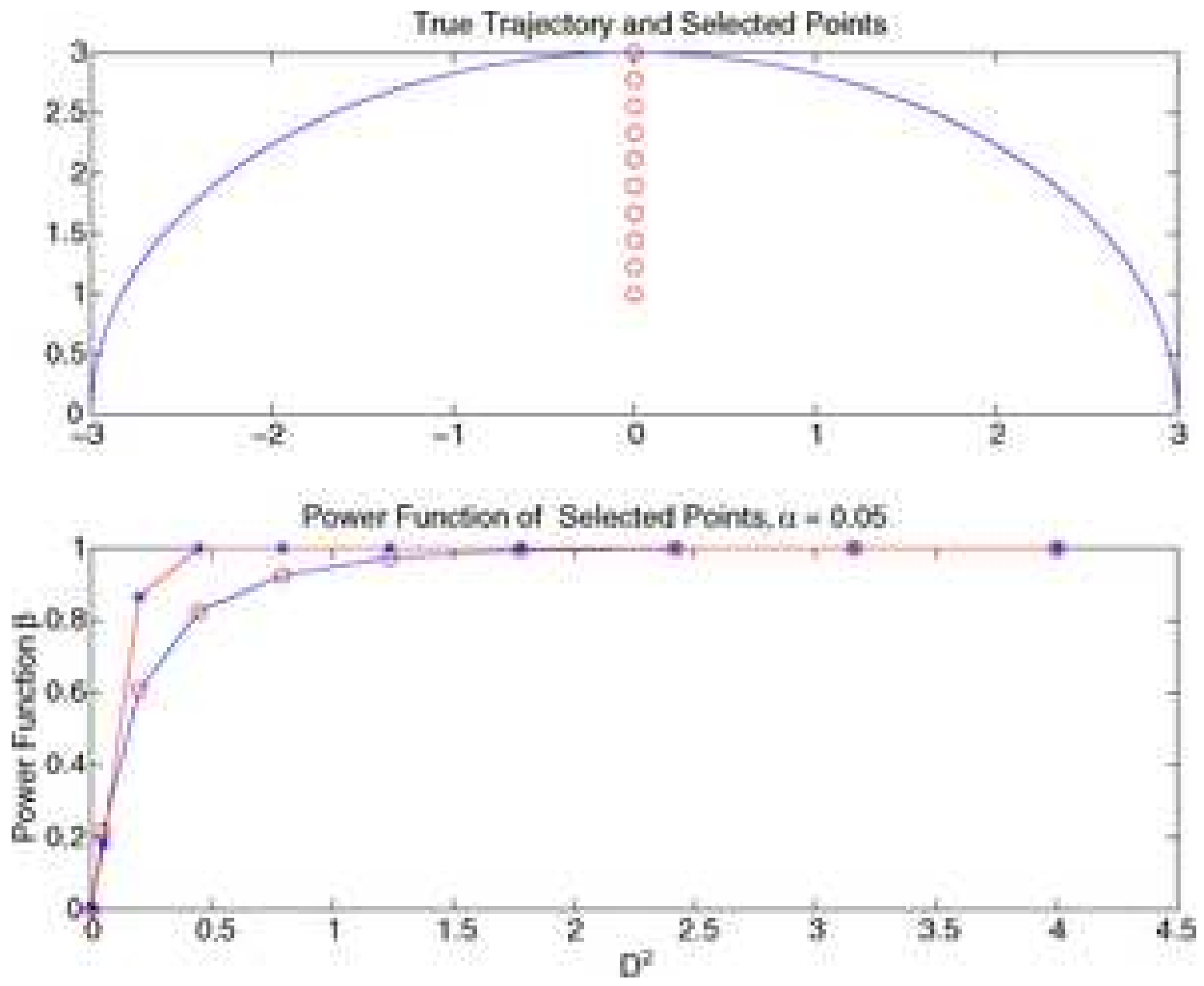

FIG. 6. *The top part shows the true integral curve and selected points of interest for measuring the power. The bottom part represents graphs of the power function based on Monte Carlo study (the red curve) and based on the theoretical formula (4.7) (the blue curve).*

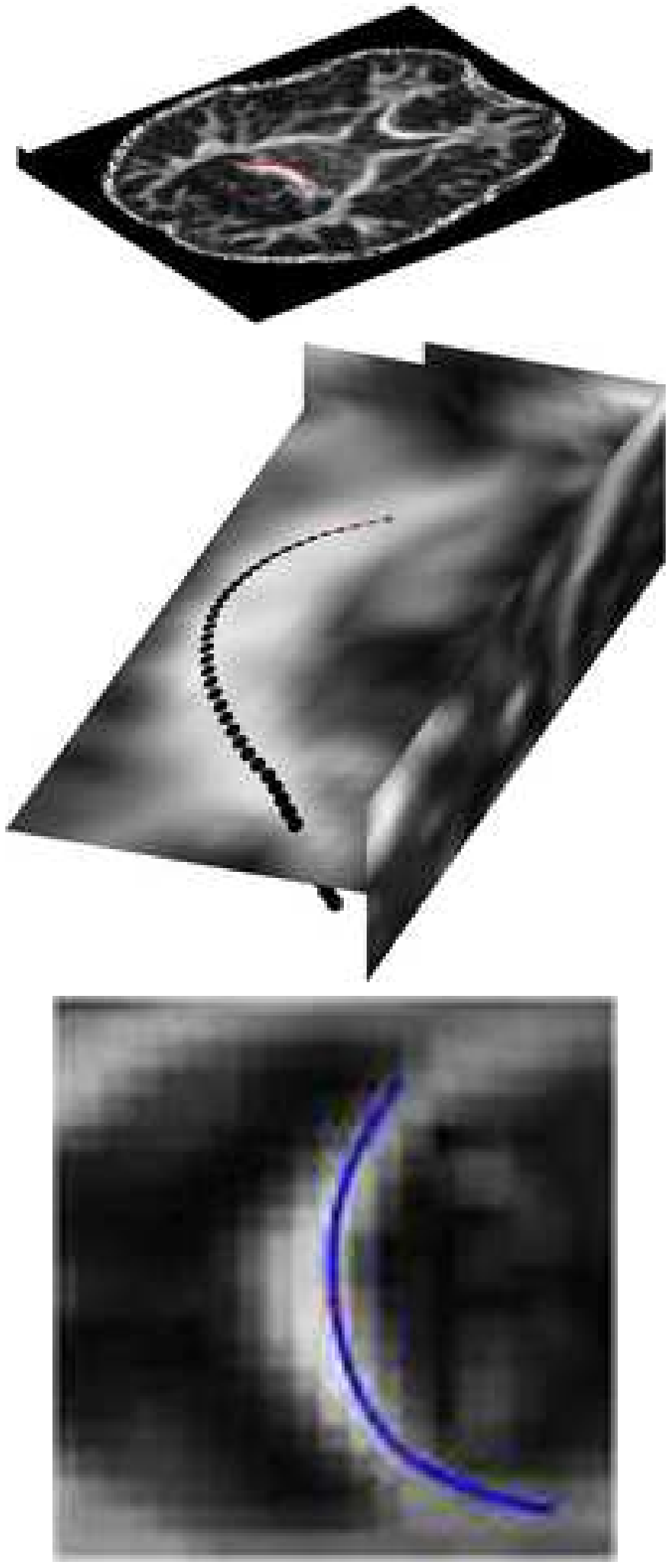

(a)

(b)

(c)

Fig. 7. *Illustrates the use of the proposed tracking procedure on real DTI data.* (a) *presents a single estimated fiber trajectory. The blue point shows the starting seed point.* (b) *shows the visualization of the 3-D confidence ellipsoid (C.E.) of the tracking procedure.* (c) *is an enlarged visualization of a 3-D fiber tracking trajectory in a given region of the brain.*

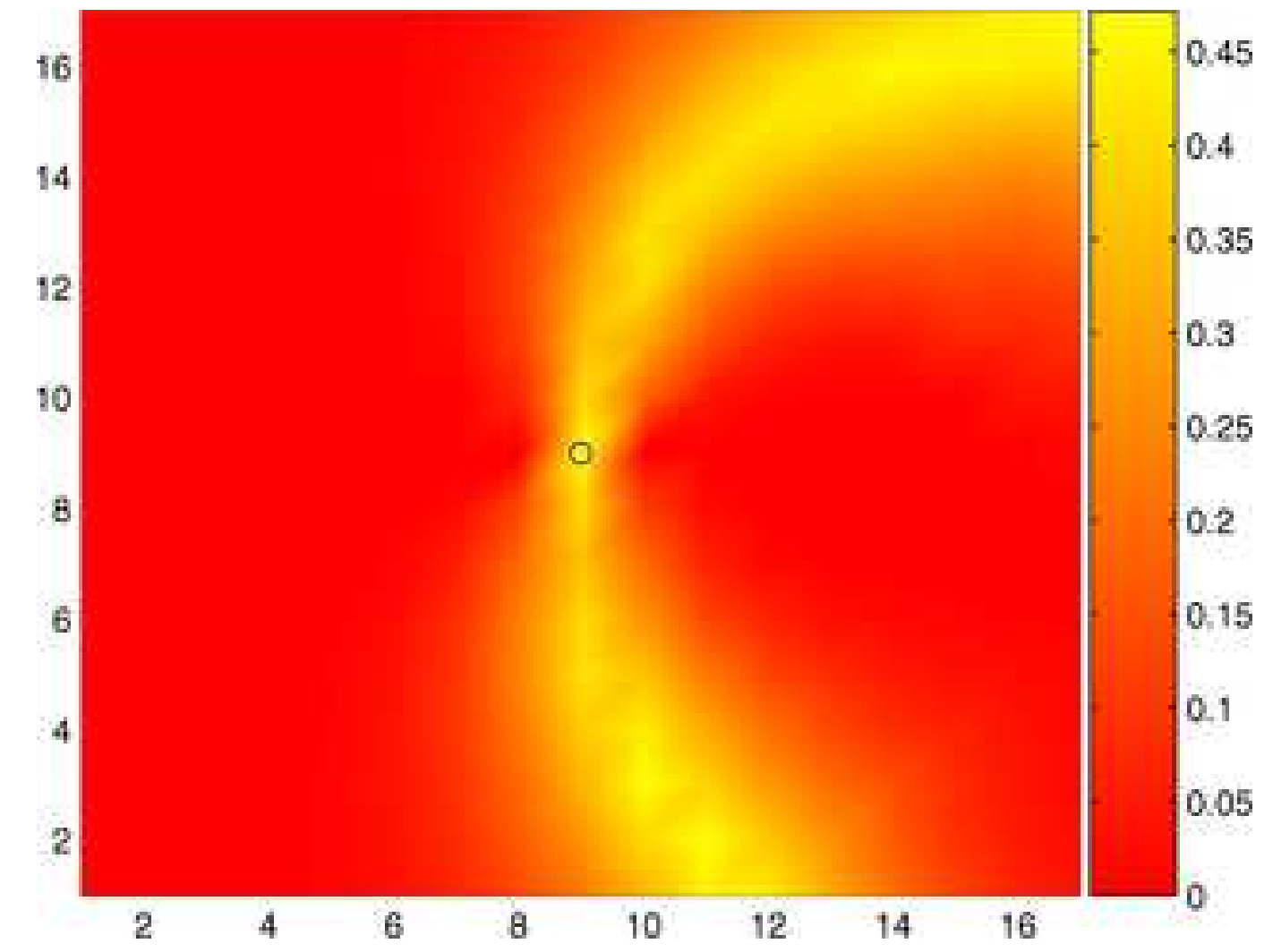

FIG. 8. *p-value map of the trajectory shown in Figure 7(c): for each point it shows the p-value of testing the null hypothesis that the true integral curve passes through this point.*

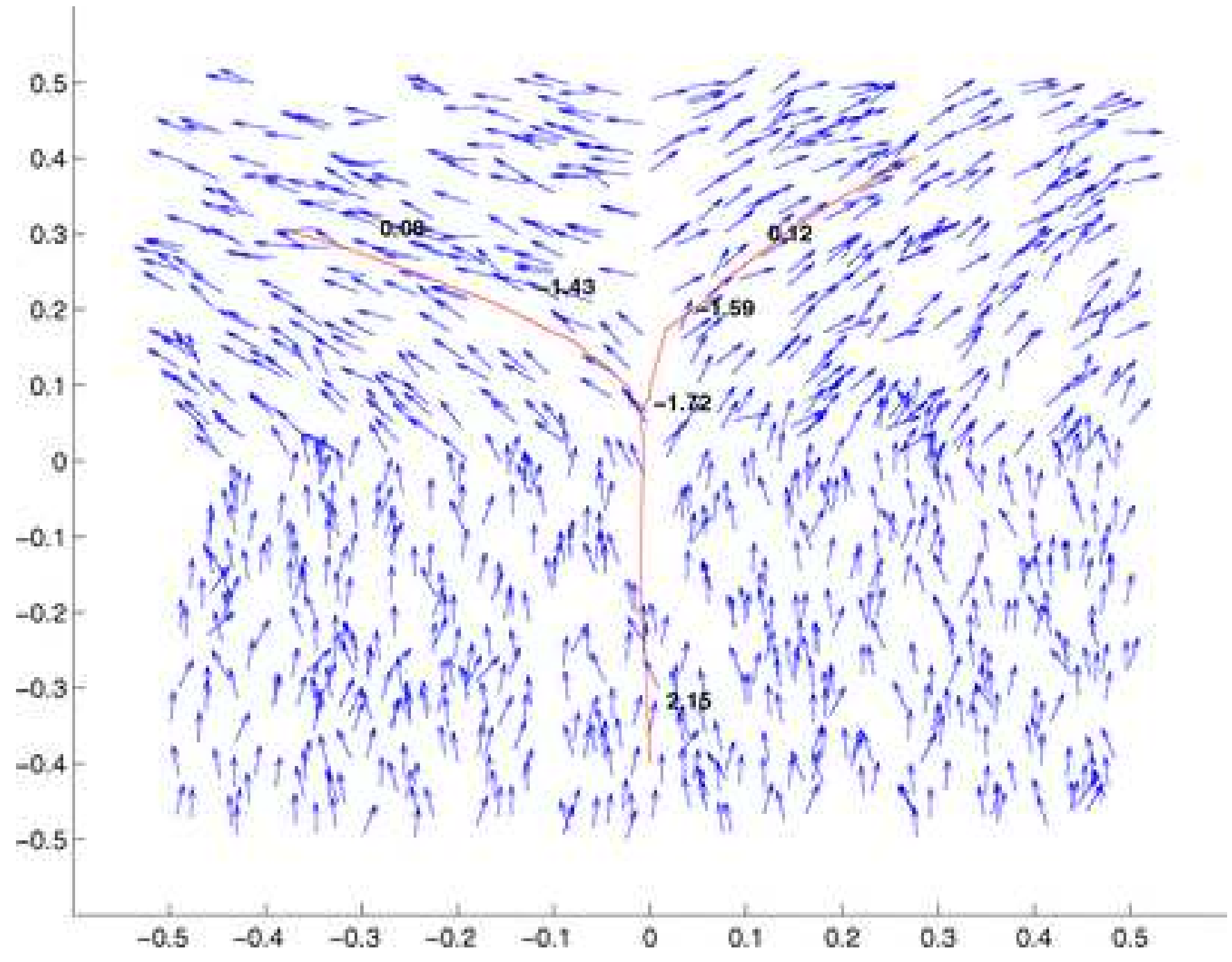

FIG. 9. *Shows the results of tracking a branching trajectory (for simulated 2-D data).*

becomes asymptotically standard normal and the test can be designed in the usual way. However, it should be pointed out that the convergence rate in this asymptotic result is rather slow and the impact of the remainder terms is significant. Moreover, when instead of the additive noise model one considers a more realistic model in which $|V_i| = 1$ (this is the case in DTI applications), the asymptotic distribution of $\nu_n$ is no longer standard normal.

The computation of the test statistic can be included in the fiber tracking algorithm. If the test detects an intersection at some moment $t$, the computation of the weighted local average $\hat{V}(\hat{X}(t))$ can be replaced by local clustering of vectors $V_i$ around this point. For instance, one can use for this purpose a weighted version of $k$-means clustering. As a result, there will be more than one continuation of the trajectory from the point $\hat{X}(t)$. A simulated example of testing and tracking a branching trajectory using this approach is given in Figure 9. The branching point of the true trajectory is $(0,0)$. The figure also shows the values of the test statistic at several points along the trajectory.

A detailed discussion of this and other applications of our methodology to DTI and its comparison with other methods goes beyond the scope of this paper and will be given in further publications in more specialized journals on neuroimaging.

**Acknowledgments.** The authors are thankful to the referees for several interesting suggestions.

V. KOLTCHINSKII
SCHOOL OF MATHEMATICS
GEORGIA INSTITUTE OF TECHNOLOGY
ATLANTA, GEORGIA 30332-0160
USA
E-MAIL: vlad@math.gatech.edu

L. SAKHANENKO
DEPARTMENT OF STATISTICS AND PROBABILITY
MICHIGAN STATE UNIVERSITY
EAST LANSING, MICHIGAN 48824-1027
USA
E-MAIL: luda@stt.msu.edu

S. Cai
The Mind Imaging Center
The University of New Mexico
Albuquerque, New Mexico 87131
USA
E-mail: shcai@unm.edu